\documentclass[11pt,a4paper]{article}
\usepackage{a4wide}
\usepackage{amsmath,amssymb}
\usepackage[english]{babel}
\usepackage{graphicx}
\usepackage[numbers]{natbib}
\usepackage{theorem}
\usepackage{bm}
\usepackage{titlesec}
\usepackage{ifpdf}
\usepackage{hyperref}
\usepackage{setspace}

\hypersetup{
    bookmarksopen=false,
    bookmarksnumbered=true,
    pdftitle={A Polling Model with Smart Customers},
    pdfauthor={M.A.A. Boon, A.C.C. van Wijk, I.J.B.F. Adan, O.J. Boxma}
}
\ifpdf
  \hypersetup{colorlinks=true,linkcolor=black,urlcolor=black,citecolor=black,pdfpagemode=UseOutlines,plainpages=false,pdfpagelabels}
\else
  \hypersetup{colorlinks=false}
\fi

\titleformat{\subsubsection}[hang]{\normalfont\em}{}{1ex}{}

\theorembodyfont{\upshape}
\theoremheaderfont{\bfseries}

\newtheorem{theorem}{Theorem}[section]
\newtheorem{property}[theorem]{Property}
\newtheorem{lemma}[theorem]{Lemma}
\newtheorem{remark}[theorem]{Remark}

\newtheorem{corollary}[theorem]{Corollary}
\newenvironment{proof}{\par\textbf{Proof\ }\ }{\hfill$\square$\par}
\numberwithin{equation}{section}

\newcommand{\ee}{\textrm{e}}

\newcommand{\z}{\mathbf{z}}
\newcommand{\zi}{\mathbf{z_i}}
\newcommand{\E}{\mathbb{E}}
\renewcommand{\P}{\mathbb{P}}

\newcommand{\BP}{\textit{BP}}
\providecommand{\href}[2]{#2}
\newcommand{\LB}{\widetilde{\textit{LB}}}
\newcommand{\VB}{\widetilde{\textit{VB}}}
\newcommand{\VC}{\widetilde{\textit{VC}}}

\newcommand{\LQ}{\hat{L}}
\newcommand{\EintForw}[2]{\E[\overrightarrow{#1}^{(#2)}]}
\newcommand{\EintBack}[2]{\E[\overleftarrow{#1}^{(#2)}]}
\newcommand{\equaldist}{\,{\buildrel d \over =}\,}
\newcommand{\indicator}{\mathbf{1}_{[\textit{LB}_i\begin{picture}(0,0)\put(-0.9,0){$\mbox{}^{(V_{j})}$}\end{picture}\phantom{(V_{j})}=\,k]}}

\title{A Polling Model with Smart Customers\footnote{The research was done in the framework of the BSIK/BRICKS project, and of the European Network of Excellence Euro-NF.}}
\author{M.A.A. Boon\footnote{\textsc{Eurandom} and Department of Mathematics and Computer Science, Eindhoven University of Technology, P.O. Box 513, 5600MB Eindhoven, The Netherlands}\\\href{mailto:marko@win.tue.nl}{marko@win.tue.nl} \and A.C.C. van Wijk\footnote{\textsc{Eurandom}, Department of Industrial Engineering \& Innovation Sciences and Department of Mathematics and Computer Science, Eindhoven University of Technology, P.O. Box 513, 5600MB Eindhoven, The Netherlands} \\\href{mailto:a.c.c.v.wijk@tue.nl}{a.c.c.v.wijk@tue.nl}\and I.J.B.F. Adan\footnotemark[2]\\\href{mailto:iadan@win.tue.nl}{iadan@win.tue.nl} \and O.J. Boxma\footnotemark[2]\\\href{mailto:boxma@win.tue.nl}{boxma@win.tue.nl}}

\date{September, 2010}

\begin{document}
\maketitle

\begin{abstract}
In this paper we consider a single-server, cyclic polling system with switch-over times. A distinguishing feature of the model is that the rates of the Poisson arrival processes at the various queues depend on the server location. For this model we study the joint queue length distribution at polling epochs and at server's departure epochs. We also study the marginal queue length distribution at arrival epochs, as well as at arbitrary epochs (which is not the same in general, since we cannot use the PASTA property). A generalised version of the distributional form of Little's law is applied to the joint queue length distribution at customer's departure epochs in order to find the waiting time distribution for each customer type. We also provide an alternative, more efficient way to determine the mean queue lengths and mean waiting times, using Mean Value Analysis. Furthermore, we show that under certain conditions a Pseudo-Conservation Law for the total amount of work in the system holds. Finally, typical features of the model under consideration are demonstrated in several numerical examples.

\end{abstract}

\section{Introduction}\label{introduction}

The classical polling system is a queueing system consisting of multiple queues, visited by a single server.
Typically, queues are served in cyclic order, and switching from one queue to the next queue requires a switch-over time, but these assumptions are not essential to the analysis. The decision at what moment the server should start switching to the next queue is important to the analysis, though. Polling systems satisfying a so-called \emph{branching property} generally allow for an exact analysis, whereas polling systems that do not satisfy this property rarely can be analysed in an exact way. See Resing \cite{resing93}, or Fuhrmann \cite{fuhrmann81}, for more details on this branching property.

There is a huge literature on polling systems, mainly because of their practical relevance. Applications are found, among others, in production environments, transportation, and data communication. The surveys of Takagi \cite{takagi1988qap}, Levy and Sidi \cite{levysidi90}, and Vishnevskii and Semenova \cite{vishnevskiisemenova06} provide a good overview of applications of polling systems. These surveys, and \cite{winandsPhD}, Chapters $2.2$ and $3$, are also excellent references to find more information about various analysis techniques, such as the Buffer Occupancy method, the Descendant Set approach, and Mean Value Analysis (MVA) for polling systems. The vast majority of papers on polling models assumes that the arrival rate stays constant throughout a cycle, although it may vary per queue. The polling model considered in the present paper, allows the arrival rate in each queue to vary depending on the server location. This model was first considered by Boxma \cite{smartcustomers}, who refers to this model as a polling model with \emph{smart customers}, because one way to look at this system is to regard it as a queueing system where customers choose which queue to join, based on the current server position. Note that Boxma's definition of smart customers is different from the definition used by Mandelbaum and Yechiali \cite{mandelbaumyechiali83}, who study an $M/G/1$ queue where smart customers may decide upon arrival to join the queue, not to enter the system at all, or to wait for a while and postpone the decision.
%Allowing arrival rates to depend on the location of the server has practical relevance because in, e.g., certain production environments or traffic intersections the arrival rates are influenced by the position of the server.

\begin{figure}[ht]
\begin{center}
\includegraphics[width=0.5\linewidth]{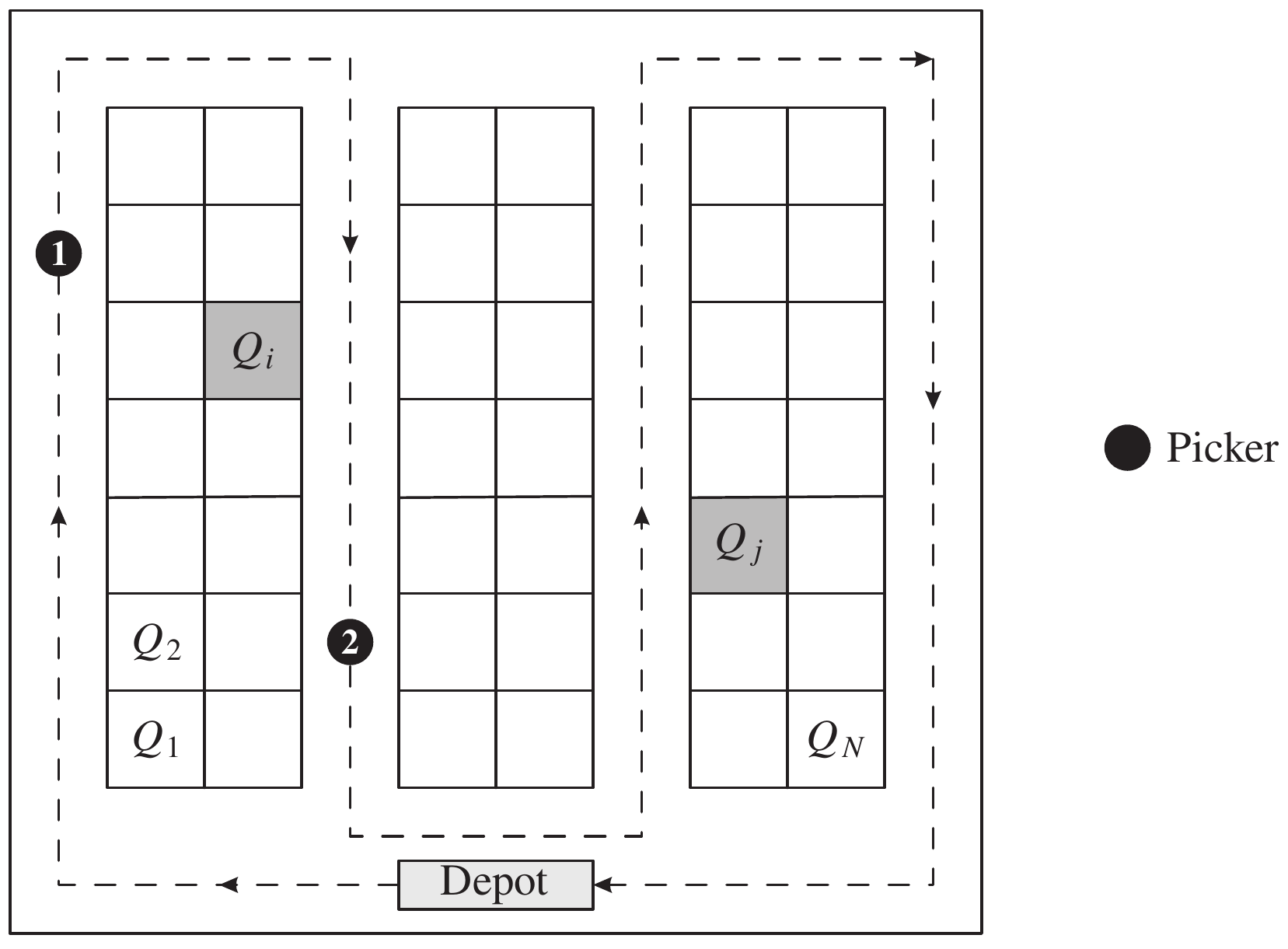}
\end{center}
\caption{A dynamic order picking system. Orders are placed in queues $Q_1,\dots,Q_N$.}
\label{zonepicking}
\end{figure}

A relevant application can be found in \cite{gongdekoster08}, where a polling model is used to model a dynamic order picking system (DPS). In a DPS, a worker picks orders arriving in real time during the picking operations and the picking information can dynamically change in a picking cycle.
One of the challenging questions that online retailers now face, is how to organise the logistic fulfillment processes during and after order receipt. In traditional stores, purchased products can be taken home immediately. However, in the case of online retailers, the customer must wait for the shipment to arrive. In order to reduce throughput times, an efficient enhancement to an ordinary DPS is to have products stored at multiple locations. The system can be modelled as a polling system with queues corresponding to each of the locations, and customers corresponding to orders. The location of the worker determines in which of the queues an order is being placed. In this system arrival rates of the orders depend on the location of the server (i.e. the worker), which makes it a typical smart customers example.
A graphical illustration is given in Figure \ref{zonepicking}. We focus on one specific order type, which is placed in two locations, say $Q_i$ and $Q_j$. While the picker is on its way to $Q_i$, say at location 1, all of these orders are routed to $Q_i$ and the arrival rate at $Q_j$ is zero. If the picker is between $Q_i$ and $Q_j$, say at location 2, the situation is reversed and $Q_j$ receives all of these orders.

Besides practical relevance, the smart customers model also provides a powerful framework to analyse more complicated polling models. For example, a polling model where the service discipline switches each cycle between gated and exhaustive, can be analysed constructing an alternative polling model with twice the number of queues and arrival rates being zero during specific visit periods \cite{boxmavanwijkadan2008}. The idea of temporarily setting an arrival rate to zero is also used in \cite{boonadan2009} for the analysis of a polling model with multiple priority levels. Time varying arrival rates also play a role in the analysis of a polling model with reneging at polling instants \cite{boonreneging2009}.

Concerning state dependent arrival rates, more literature is available for systems consisting of only one queue, often assuming phase-type distributions for vacations and/or service times. A system consisting of a single queue with server breakdowns and arrival rates depending on the server status is studied in \cite{shogan79}. A difference with the system studied in the present paper, besides the number of queues, is that the machine can break down at arbitrary moments during the service of customers. Polling systems with breakdowns have been studied as well, cf. \cite{serverbreakdowns1,stationandserverbreakdowns,stationbreakdowns1,stationbreakdowns2}. However, only Nakdimon and Yechiali \cite{stationbreakdowns2} consider a model where the arrival process stops temporarily during a breakdown.
Shanthikumar \cite{shanthikumar88} discusses a stochastic decomposition for the queue length in an $M/G/1$ queue with server vacations under less restrictive assumptions than Fuhrmann and Cooper \cite{fuhrmanncooper85}. One of the relaxations is that the arrival rate of customers may be different during visit periods and vacations.
Another system, with so-called working vacations and server breakdowns is studied in \cite{jainjain2010}. During these working vacations, both the service and arrival rates are different. Mean waiting times are found using a matrix analytical approach.
For polling systems, a model with arrival rates that vary depending on the location of the server has not been studied in detail yet. Boxma \cite{smartcustomers} studies the joint queue length distribution at the beginning of a cycle, but no waiting times or marginal queue lengths are discussed. In a recent paper \cite{pollinglevy09}, a polling system with L\'evy-driven, possibly correlated input is considered. Just as in the present paper, the arrival process may depend on the location of the server. In \cite{pollinglevy09} typical performance measures for L\'evy processes are determined, such as the steady-state distribution of the joint amount of fluid at an arbitrary epoch, and at polling and switching instants.
The present paper studies a similar setting, but assumes Poisson arrivals of individual customers. This enables us to find the probability generating functions (PGFs) of the joint queue length distributions at polling instants and customer's departure epochs, and the marginal queue length distributions at customer's arrival epochs and at arbitrary epochs (which are not the same, because PASTA cannot be used).
The introduction of customer subtypes, categorised by their moment of arrival, makes it possible to generalise the distributional form of Little's law (see, e.g., \cite{keilsonservi90}), and apply it to the joint queue length distribution at departure epochs to find the Laplace-Stieltjes Transform (LST) of the waiting time distribution.

The present paper is structured as follows: Section \ref{notation} gives a detailed model description and introduces the notation used in this paper. In Section \ref{nosubtypes} the PGFs of the joint queue length distributions of all customer types at polling instants are derived. The marginal queue length distribution is also studied in this section, but we show in Section \ref{subtypes} that the derivation of the waiting time LST for each customer type requires a more complicated analysis, based on customer subtypes. In Sections \ref{nosubtypes} and \ref{subtypes} we need information on the lengths of the cycle time and all visit times, which are studied in Section~\ref{cycletimesection}. In Section~\ref{mva} we adapt the MVA framework for polling systems, introduced in \cite{winands06}, to our model. This results in a very efficient method to compute the mean waiting time of each customer type. For polling systems with constant arrival rates, a Pseudo-Conservation Law (PCL) is studied in \cite{boxmagroenendijk87}. In Section \ref{pclsection} we show that, under certain conditions, a PCL is satisfied by our model. Finally, we give numerical examples that illustrate some typical features and advantages of the model under consideration.

\section{Model description and notation\label{notation}}

The polling model in the present paper contains $N$ queues, $Q_1, \dots, Q_N$, visited in cyclic order by one server. Switching from $Q_i$ to $Q_{i+1}$ ($i=1,\dots,N$, where $Q_{N+1}$ is understood to be $Q_1$, etc.) requires a switch-over time $S_i$, with LST $\sigma_i(\cdot)$. We assume that at least one switch-over time is strictly greater than zero, otherwise the mean cycle length in steady-state becomes zero and the analysis changes slightly. See, e.g., \cite{borst97} for a relation between polling systems with and without switch-over times. Switch-over times are assumed to be independent. The cycle time $C_i$ is the time that elapses between two successive visit \emph{beginnings} to $Q_i$, and $C^*_i$ is the time that elapses between two successive visit \emph{endings} to $Q_i$. The mean cycle time does not depend on the starting point of the cycle, so $\E[C_i] = \E[C^*_i] =\E[C]$. The visit time $V_i$ of $Q_i$ is the time between the visit beginning and visit ending of $Q_i$.
The intervisit time $I_i$ of $Q_i$ is the time between a visit \emph{ending} to $Q_i$ and the next visit \emph{beginning} at $Q_i$. We have $C_i = V_i + I_i$, and $I_i = S_i + V_{i+1} + \dots + S_{i+N-1}$, $i=1,\dots,N$. Customers arriving at $Q_i$, i.e. type $i$ customers, have a service requirement $B_i$, with LST $\beta_i(\cdot)$. We also assume independence of service times, and first-come-first-served (FCFS) service order.

The service discipline of each queue determines the moment at which the server switches to the next queue. In the present paper we study the two most popular service disciplines in polling models, exhaustive service (the server switches to the next queue directly after the last customer in the current queue has been served) and gated service (only visitors present at the server's arrival at the queue are served). The reason why these two service disciplines have become the most popular in polling literature, lies in the fact that they are from a practical point of view the most relevant service disciplines that allow an exact analysis. In this respect the following property, defined by Resing \cite{resing93} and also Fuhrmann \cite{fuhrmann81}, is very important.
\begin{property}[Branching Property]\label{resingproperty}
If the server arrives at $Q_i$ to find $k_i$ customers there, then during the course of the server's visit, each of these $k_i$ customers will effectively be replaced in an i.i.d. manner by a random population having probability generating function $h_i(z_1,\dots,z_N)$, which can be any $N$-dimensional probability generating function.
\end{property}
In most cases, a polling model can only be analysed exactly, if the service discipline at each queue satisfies Property \ref{resingproperty}, or some slightly weaker variant of this property, because in this case the joint queue length process at visit beginnings to a fixed queue constitutes a Multi-Type Branching Process, which is a nicely structured and well-understood process. Gated and exhaustive service both satisfy this property, whereas a service discipline like $k$-limited service (serve at most $k$ customers during each visit) does not.

The feature that distinguishes the model under consideration from commonly studied polling models, is the arrival process.
This arrival process is a standard Poisson process, but the rate depends on the location of the server. The arrival rate at $Q_i$ is denoted by $\lambda_i^{(P)}$, where $P$ denotes the position of the server, which is either serving a queue, or switching from one queue to the next: $P \in \{V_1, S_1, \dots, V_N, S_N\}$. One of the consequences is that the PASTA property does not hold for an arbitrary arrival, but as we show in Section \ref{nosubtypes}, a conditional version of PASTA does hold. Another difficulty that arises, is that the distributional form of Little's law cannot be applied to the PGF of the marginal queue length distribution to obtain the LST of the waiting time distribution anymore. We explain this in Section \ref{subtypes}, where we also derive a generalisation of the distributional form of Little's law.

\section{Queue length distributions\label{nosubtypes}}

\subsection{Joint queue length distribution at visit beginnings/endings}

The two main performance measures of interest, are the steady-state queue length distribution and the waiting time distribution of each customer type. In this section we focus on queue lengths rather than waiting times, because the latter requires a more complex approach that is discussed in the next section. We restrict ourselves to branching-type service disciplines, i.e., service disciplines satisfying Property \ref{resingproperty}. Boxma \cite{smartcustomers} follows the approach by Resing \cite{resing93}, defining offspring and immigration PGFs to determine the joint queue length distribution at the beginning of a cycle. We take a slightly different approach that gives the same result, but has the advantage that it gives expressions for the joint queue length PGF at all visit beginnings and endings as well. %Denote by $V_{b_i}(z_1, \dots,z_N)$ the PGF of the steady-state joint queue length distribution at visit \emph{beginnings} to $Q_i$. Similarly, $V_{c_i}(z_1, \dots,z_N)$ is the equivalent at visit \emph{endings}.
Denote by $\LB^{(P)}(z_1, \dots, z_N)$ the PGF of the steady-state joint queue length distribution at beginnings of period $P \in \{V_1, S_1, \dots, V_N, S_N\}$.
The relation between these PGFs, also referred to as \emph{laws of motion} in the polling literature, is obtained by application of Property \ref{resingproperty} to $\LB^{(V_i)}(\z)$, where $\z$ is a shorthand notation for the vector $(z_1,\dots, z_{N})$. This property states that each type $i$ customer present at the visit beginning to $Q_i$ will be replaced during this visit by a random population having PGF $h_i(\z)$, which depends on the service discipline. The only difference between conventional polling models and the model under consideration in the present paper, is that the arrival rates depend on the server location. The relations between $\LB^{(V_i)}(\z), \LB^{(S_i)}(\z)$, and $\LB^{(V_{i+1})}(\z)$ are given by:
\begin{align}
\LB^{(S_i)}(\z) &= \LB^{(V_i)}(z_1, \dots,z_{i-1}, h_i(\z), z_{i+1}, \dots, z_N),\label{VCi}\\
\LB^{(V_{i+1})}(\z) &= \LB^{(S_i)}(\z)\,\sigma_i\Big(\sum_{j=1}^N \lambda_j^{(S_i)}(1-z_j)\Big),\label{VBi}
\end{align}
where $h_{i}(\z)$ is the PGF mentioned in Property \ref{resingproperty}. It is discussed in the context of a polling model with smart customers in \cite{smartcustomers}. For gated service, $h_i(\z) = \beta_i\left(\sum_{j=1}^N\lambda_j^{(V_i)}(1-z_j)\right)$. For exhaustive service, $h_i(\z) = \pi_i\left(\sum_{j\neq i}\lambda_j^{(V_i)}(1-z_j)\right)$, where $\pi_i(\cdot)$ is the LST of a busy period distribution in an $M/G/1$ system with only type $i$ customers, so it is the root in $(0,1]$ of the equation $\pi_i(\omega) = \beta_i\left(\omega + \lambda_i^{(V_i)}(1 - \pi_i(\omega))\right)$, $\omega \geq 0$ (cf. \cite{cohen82}, p. 250). Now that we can relate $\LB^{(V_{i+1})}(\z)$ to $\LB^{(V_{i})}(\z)$, we can repeat this and finally obtain a recursion for $\LB^{(V_{i})}(\z)$. This recursive expression is sufficient to compute all moments of the joint queue length distribution at a visit beginning to $Q_{i}$ by differentiation, but iteration of the expression leads to the steady-state queue length distribution at polling epochs, written as an infinite product. We refer to \cite{resing93} for more details regarding this approach, and for rigorous proofs of the laws of motion. Stability conditions are studied in more detail in \cite{pollinglevy09}, where it is shown that a necessary and sufficient condition for ergodicity is that
the Perron-Frobenius eigenvalue of the matrix $R - I_N$ should be less than 0, where $I_N$ is the $N \times N$ identity matrix, and $R$ is an $N \times N$ matrix containing elements $\rho_{ij} := \lambda_{i}^{(V_j)}\E[B_i]$. This holds under the assumption that $\E[V_i] > 0$ for all $i=1,\dots,N$.

\subsection{Marginal queue length distribution}

Common techniques in polling systems (see, e.g. \cite{semphd,eisenberg72}) to determine the PGF of the steady-state marginal queue length distribution of each customer type, are based on deriving the queue length distribution at departure epochs. A level-crossing argument implies that the marginal queue length distribution at arrival epochs must be the same as the one at departure epochs, and, finally, because of PASTA this distribution is the same as the marginal queue length distribution at an arbitrary point in time. In our model, the marginal queue length distributions at arrival and departure epochs are also the same, but the distribution at arbitrary moments is different because of the varying arrival rates during a cycle. We can circumvent this problem by conditioning on the location $P$ of the server $(P \in \{V_1, S_1, \dots, V_N, S_N\})$ and use conditional PASTA to find the PGF of the marginal queue length distribution at an arbitrary point in time. Let $L_i$ denote the steady-state queue length of type $i$ customers at an arbitrary moment, and let $L_i^{(V_j)}$ and $L_i^{(S_j)}$ denote the queue length of type $i$ customers at an arbitrary time point during $V_j$ and $S_j$ respectively $(i,j = 1,\dots,N)$. The following relation holds:
\begin{equation}
\E[z^{L_i}] = \sum_{j=1}^N\left(\frac{\E[V_j]}{\E[C]}\E\left[z^{L_i^{(V_j)}}\right] + \frac{\E[S_j]}{\E[C]}\E\left[z^{L_i^{(S_j)}}\right]\right), \qquad i=1,\dots,N.
\label{queuelengthGF}
\end{equation}
Note that, at this moment, $\E[V_j]$ and $\E[C]$ are still unknown. In Sections \ref{cycletimesection} and \ref{mva} we illustrate two different ways to compute them.
Since $S_j$, for $j=1,\dots,N$, and $V_j$, for $j\neq i$, are \emph{non-serving} intervals for customers of type $i$, we use a standard result (see, e.g., \cite{semphd}) to find the PGFs of $L_i^{(V_j)}$ and $L_i^{(S_j)}$ respectively:
\begin{align}
\E\left[z^{L_i^{(V_j)}}\right] &= \frac{\E[z^{\textit{LB}_i^{(V_{j})}}]-\E[z^{\textit{LB}_i^{(S_{j})}}]}{(1-z)\left(\E[\textit{LB}_i^{(S_{j})}]-\E[\textit{LB}_i^{(V_{j})}]\right)},\qquad&&i=1,\dots,N; j\neq i,\label{queuelengthGFduringvisitj}\\
\E\left[z^{L_i^{(S_j)}}\right] &= \frac{\E[z^{\textit{LB}_i^{(S_{j})}}]-\E[z^{\textit{LB}_i^{(V_{j+1})}}]}{(1-z)\left(\E[\textit{LB}_i^{(V_{j+1})}]-\E[\textit{LB}_i^{(S_{j})}]\right)},\qquad&&i,j=1,\dots,N,\label{queuelengthGFduringswitchoverj}
\end{align}
%where $\textit{LB}_i^{(V_{j})}$ and $\textit{LB}_i^{(S_{j})}$ are the number of type $i$ customers at respectively a visit beginning and ending at $Q_j$.
where $\textit{LB}_i^{(P)}$, for $i=1,\dots,N$, are the number of type $i$ customers at the beginning of period $P \in \{V_1, S_1, \dots, V_N, S_N\}$.
Their PGFs can be expressed in terms of $\LB^{(V_{1})}(\z)$ using the relations \eqref{VBi} and \eqref{VCi}, and replacing argument $\z$ by the vector $(1,\dots,1,z,1,\dots,1)$ where $z$ is the element at position $i$. Using branching theory from \cite{resing93}, an explicit expression for $\LB^{(V_{1})}(\z)$ is given in \cite{smartcustomers}.
The mean values, $\E[\textit{LB}_i^{(V_{j})}]$ and $\E[\textit{LB}_i^{(S_{j})}]$, can be obtained by differentiation of the corresponding PGFs and substituting $z = 1$. %To compute these means, no explicit computation of $V_{b_1}(\z)$ is required.

It remains to compute $\E\left[z^{L_i^{(V_i)}}\right]$, $i=1,\dots,N$, i.e., the PGF of the number of type $i$ customers at an arbitrary point within $V_i$. As far as the marginal queue length of type $i$ customers is concerned, the system can be viewed as a vacation queue with the intervisit time $I_i$ corresponding to the server vacation. We can use the Fuhrmann-Cooper decomposition \cite{fuhrmanncooper85}, but we have to be careful here. In a polling system where type $i$ customers arrive with \emph{constant arrival rate} $\lambda_i^{(V_i)}$, the Fuhrmann-Cooper decomposition states that
\begin{equation}
\E[z^{L_i}] = \frac{(1-\lambda_i^{(V_i)}\E[B_i])(1-z)\beta_i\big(\lambda_i^{(V_i)}(1-z)\big)}{\beta_i\big(\lambda_i^{(V_i)}(1-z)\big)-z} \times \frac{\E\left[z^{\textit{LB}_i^{(S_{i})}}\right]-\E\left[z^{\textit{LB}_i^{(V_{i})}}\right]}{(1-z)\left(\E[\textit{LB}_i^{(V_{i})}]-\E[\textit{LB}_i^{(S_{i})}]\right)}.\label{fuhrmanncooperdecomposition}
\end{equation}
The two parts in this decomposition can be recognised as the PGFs of the number of type $i$ customers respectively at an arbitrary moment in an $M/G/1$ queue, and at an arbitrary point during the intervisit time $I_i$. Of course, the following relation also holds:
\begin{equation}
\E[z^{L_i}] = \frac{\E[V_i]}{\E[C]}\E[z^{L_i^{(V_i)}}] + \frac{\E[I_i]}{\E[C]}\E[z^{L_i^{(I_i)}}].\label{EzLi2}
\end{equation}
Combining \eqref{fuhrmanncooperdecomposition} with \eqref{EzLi2}, results in:
\begin{equation}
\E[z^{L_i^{(V_i)}}] = \frac{1-\lambda_i^{(V_i)}\E[B_i]}{\lambda_i^{(V_i)}\E[B_i]} \frac{z\big(1-\beta_i(\lambda_i^{(V_i)}(1-z))\big)}{\beta_i(\lambda_i^{(V_i)}(1-z))-z}
\times\frac{\E\left[z^{\textit{LB}_i^{(S_{i})}}\right]-\E\left[z^{\textit{LB}_i^{(V_{i})}}\right]}{(1-z)\left(\E[\textit{LB}_i^{(V_{i})}]-\E[\textit{LB}_i^{(S_{i})}]\right)},\label{queuelengthGFduringvisiti}
\end{equation}
for $i=1,\dots,N$. The second part of this decomposition is, again, the PGF of the number of customers at an arbitrary point during the intervisit time $I_i$. The first part can be recognised as the PGF of the queue length of an $M/G/1$ queue with type $i$ customers at an arbitrary point \emph{during} a busy period.

Now we return to the model \emph{with} varying arrival rates. The key observation is that the behaviour of the number of type $i$ customers \emph{during a visit period} of $Q_i$, is exactly the same in this system as in a polling system with constant arrival rates $\lambda_i^{(V_i)}$ for type $i$ customers. Equation \eqref{queuelengthGFduringvisiti} no longer depends on anything that happens during the intervisit time, because this is all captured in $\textit{LB}_i^{(V_{i})}$, the number of type $i$ customers at the beginning of a visit to $Q_i$. This implies that, for a polling model with smart customers, the queue length PGF of $Q_i$ at a random point during $V_i$ is also given by \eqref{queuelengthGFduringvisiti}. The only difference lies in the interpretation of \eqref{queuelengthGFduringvisiti}. Obviously, the first part in \eqref{queuelengthGFduringvisiti} is still the PGF of the queue length distribution of an $M/G/1$ queue at an arbitrary point during a busy period. However, the last term can no longer be interpreted as the PGF of the distribution of the number of type $i$ customers at an arbitrary point during the intervisit time $I_i$.

Substitution of \eqref{queuelengthGFduringvisitj}, \eqref{queuelengthGFduringswitchoverj}, and \eqref{queuelengthGFduringvisiti} in \eqref{queuelengthGF} gives the desired expression for the PGF of the marginal queue length in $Q_i$.

\begin{remark}
The marginal queue length PGF \eqref{queuelengthGF} has been obtained by conditioning on the position of the server at an arbitrary epoch in a cycle, which explains the probabilities $\frac{\E[V_j]}{\E[C]}$ (server is serving $Q_j$) and $\frac{\E[S_j]}{\E[C]}$ (server is switching to $Q_{j+1}$). It is easy now to obtain the marginal queue length PGF at an \emph{arrival epoch}, simply by conditioning on the position of the server at an arbitrary \emph{arrival epoch}. The probability that the server is at position $P \in \{V_1, S_1, \dots, V_N, S_N\}$ at the arrival of a type~$i$ customer, is $\frac{\lambda_i^{(P)}\E[P]}{\overline{\lambda}_i\E[C]}$, with $\overline{\lambda}_i=\frac{1}{\E[C]}\sum_{j=1}^N\left(\lambda_i^{(V_j)}\E[V_j]+\lambda_i^{(S_j)}\E[S_j]\right)$. This results in the following expression for the PGF of the distribution of the number of type $i$ customers at the arrival of a type $i$ customer:
\begin{equation}
\E[z^{L_i}|\textrm{arrival type $i$}] = \sum_{j=1}^N\left(\frac{\lambda_i^{(V_j)}\E[V_j]}{\overline{\lambda}_i\E[C]}\E\left[z^{L_i^{(V_j)}}\right] + \frac{\lambda_i^{(S_j)}\E[S_j]}{\overline{\lambda}_i\E[C]}\E\left[z^{L_i^{(S_j)}}\right]\right),
\label{queuelengthGFarrival}
\end{equation}
for  $i=1,\dots,N$. A standard up-and-down crossing argument can be used to argue that \eqref{queuelengthGFarrival} is also the PGF of the distribution of the number of type $i$ customers at the \emph{departure} of a type $i$ customer. As stated before, it is different from the PGF of the distribution of the number of type $i$ customers at an \emph{arbitrary} epoch, unless $\lambda_i^{(V_j)}=\lambda_i^{(S_j)}=\overline{\lambda}_i$ for all $i,j=1,\dots,N$ (as is the case in polling models without smart customers).
\end{remark}

\begin{remark}\label{thetadef}
Equations \eqref{queuelengthGFduringvisitj} and \eqref{queuelengthGFduringswitchoverj} rely heavily on the PASTA property and are only valid if type~$i$ arrivals take place during the non-serving interval. If no type $i$ arrivals take place (i.e. $\lambda_i^{(P)} = 0$ for the non-serving interval $P$), both the numerator and the denominator become 0. This situation has to be analysed differently. Now assume that $\lambda_i^{(P)} = 0$ for a specific customer type $i=1,\dots,N$, during a non-serving interval $P \in \{V_1, S_1, \dots, V_N, S_N\} \backslash V_i$. We now distinguish between visit periods and switch-over periods.
Let us first assume that $P$ is a switch-over time, say $S_j, j=1,\dots,N$. The length of a switch-over time is independent from the number of customers in the system, so the distribution of the number of type $i$ customers at an arbitrary point in time during $S_j$ is the same as at the beginning of $S_j$:
\[
\E\left[z^{L_i^{(S_j)}}\right] = \E\left[z^{\textit{LB}_i^{(S_{j})}}\right], \qquad i,j=1,\dots,N.
\]
The case where $P$ is a visit time, say $P=V_j$ for some $j\neq i$, requires more attention, because the length of $V_j$ depends on the number of type $j$ customers present at the visit beginning. Since this number is positively correlated with the number of customers in the other queues, we have to correct for the fact that it is more likely that a random point during an arbitrary $V_j$, falls within a long visit period (with more customers present at its beginning) than in a short visit period.
The first step, is to determine the probability that the number of type $i$ customers at an arbitrary point during $V_j$ is $k$. Since we consider the case where $\lambda_i^{(V_j)}=0$, this implies that we need the probability that the number of customers at the beginning of $V_j$ is $k$. Standard renewal arguments yield
\begin{align}
\begin{aligned}
\P[L_i^{(V_j)} = k] &= \frac{\P[\textit{LB}_i^{(V_{j})} = k]\,\E[V_j|\textit{LB}_i^{(V_{j})}=k]}{\sum_{l=0}^\infty \P[\textit{LB}_i^{(V_{j})} = l]\,\E[V_j|\textit{LB}_i^{(V_{j})}=l]}\\
&= \frac{\E[V_j\,\indicator]}{\E[V_j]},
\end{aligned}\label{renewalargument}
\end{align}
where $\mathbf{1}_{[A]}$ is the indicator function for event $A$. The first line in \eqref{renewalargument} is based on the fact that the probability is proportional to the length of visit periods $V_j$ that start with $k$ type $i$ customers, and to the number of such visit periods $V_j$. The denominator is simply a normalisation factor.

Now we can write down the expression for the number of type~$i$ customers at an arbitrary point during~$V_j$ if $\lambda_i^{(V_j)} = 0$:
\begin{align}
\E\left[z^{L_i^{(V_j)}}\right] %&= \frac{1}{\E[V_j]}\int_0^\infty \sum_{k=0}^\infty z^k t \dd \P(L_i^{(V_{b_j})}=k, V_j<t)   \nonumber\\
&= \sum_{k=0}^\infty z^k\, \P[L_i^{(V_j)} = k] \nonumber\\
&= \frac{1}{\E[V_j]} \sum_{k=0}^\infty z^k\, \E[V_j\,\indicator]\nonumber\\
&= \frac{1}{\E[V_j]} \E[V_j \sum_{k=0}^\infty z^k \,\indicator]\nonumber\\
&= \frac{1}{\E[V_j]} \E[V_j\, z^{\textit{LB}_i^{(V_{j})}}]\nonumber\\
&= \left.-\frac{1}{\E[V_j]}\frac{\partial}{\partial \omega} \E\left[z^{\textit{LB}_i^{(V_{j})}}\ee^{-\omega V_j}\right]\right|_{\omega=0},\label{queuelengthGFduringvisitjLambda0}
\end{align}
for $i=1,\dots,N$ and $j\neq i$.

Now we only need to determine $\E[z^{\textit{LB}_i^{(V_{j})}}\ee^{-\omega V_j}]$. We use the joint queue length distribution of all customers present at the beginning of $V_j$, which is given implicitly by \eqref{VBi}.
Define $\Theta_j$ as the time that the server spends at $Q_j$ due to the presence of one customer there, with LST $\theta_j(\cdot)$. For gated service $\theta_j(\cdot)=\beta_j(\cdot)$, and for exhaustive service $\theta_j(\cdot)=\pi_j(\cdot)$. The length of $V_j$, given that $l_j$ type $j$ customers are present at the visit beginning, is the sum of $l_j$ independent random variables with the same distribution as $\Theta_j$, denoted by $\Theta_{j,1}, \dots, \Theta_{j,l_j}$. The joint distribution of the number of type $i$ customers present at the beginning of $V_j$ and the length of $V_j$ is given by:
\begin{align}
\E\left[z^{\textit{LB}_i^{(V_{j})}}\ee^{-\omega V_j}\right]\, %&= \E\left[z^{L_i^{(V_{b_j})}}\ee^{-\omega (\Theta_{j,1}+\dots+\Theta_{j,L_j^{(V_{b_j})}})}\right]\\
&= \sum_{l_i=0}^\infty\sum_{l_j=0}^\infty\E\left[z^{l_i}\ee^{-\omega (\Theta_{j,1}+\dots+\Theta_{j,l_j})}\right]\P\left[\textit{LB}_i^{(V_{j})}=l_i, \textit{LB}_j^{(V_{j})}=l_j\right]\nonumber\\
&= \sum_{l_i=0}^\infty\sum_{l_j=0}^\infty z^{l_i}\E\left[\ee^{-\omega \Theta_{j,1}}\right]\times \dots \times \E\left[\ee^{-\omega \Theta_{j,l_j}}\right]\P\left[\textit{LB}_i^{(V_{j})}=l_i, \textit{LB}_j^{(V_{j})}=l_j\right]\nonumber\\
&= \sum_{l_i=0}^\infty\sum_{l_j=0}^\infty z^{l_i}\theta_j(\omega)^{l_j}\P\left[\textit{LB}_i^{(V_{j})}=l_i, \textit{LB}_j^{(V_{j})}=l_j\right]\nonumber\\
&= \LB^{(V_{j})}(1,\dots,1,z,1,\dots,1, \theta_j(\omega),1,\dots,1),\label{PGFzLiVbj}
\end{align}
where $z$ corresponds to customers in $Q_i$, and $\theta_j(\omega)$ corresponds to customers in $Q_j$. Substitution of \eqref{PGFzLiVbj} in \eqref{queuelengthGFduringvisitjLambda0} gives the desired result.

\end{remark}

\section{Waiting time distribution\label{subtypes}}

In the previous section we gave an expression for the PGF of the distribution of the steady-state queue length of a type $i$ customer at an arbitrary epoch, $L_i$. If the arrival rates do not depend on the server position, i.e. $\lambda_i^{(V_j)}=\lambda_i^{(S_j)}=\overline{\lambda}_i$ for all $i,j=1,\dots,N$, we can use the distributional form of Little's law (see, e.g., \cite{keilsonservi90}) to obtain the LST of the distribution of the waiting time of a type $i$ customer, $W_i$, $i=1,\dots,N$. Because of the varying arrival rates, there is no $\lambda_i$ for which the relation $\E[z^{L_i}] = \E\left[\ee^{-\lambda_i(1-z)(W_i+B_i)}\right]$ holds (even if we choose $\lambda_i = \overline{\lambda}_i$). In the present section, we introduce subtypes of each customer type. Each subtype is identified by the position of the server at its arrival in the system. We show that one can use a generalised version of the distributional form of Little's law that leads to the LST of the waiting time distribution of a type $i$ customer, when applied to the PGF of the \emph{joint} queue length distribution of all subtypes of a type $i$ customer. Determining this PGF requires a separate treatment of exhaustive and gated service, so results in this section do not apply to any \emph{arbitrary} branching-type service discipline.% (although the methodology can probably be applied to any branching-type service discipline).

\subsection{Joint queue length distribution at visit beginnings/endings for all subtypes\label{jointsubtypes}}

In the present section we distinguish between subtypes of type $i$ customers, arriving during different visit/switch-over periods. We define a type $i^{(P)}$ customer to be a customer arriving at $Q_i$ during $P \in \{V_1, S_1, \dots, V_N, S_N\}$. Therefore, only in this section, we define $\z$ in the following way:
\[\z = (z_1^{(V_1)},\dots,z_1^{(S_N)}, \dots,z_N^{(V_1)},\dots,z_N^{(S_N)}).\]
Note that $\z$ has $2N^2$ components: $N$ customer types times $2N$ subperiods within a cycle ($N$ visit times plus $N$ switch-over times). Let $\VB_{i}^{(P)}(\z)$ be the PGF of the joint queue length distribution of all these customer types at the moment that the server starts serving type $i$ customers that have arrived when the server was located at position~$P$. $\VC_{i}^{(P)}(\z)$ is defined equivalently for the moment that the server completes service of type $i^{(P)}$ customers.

For \emph{exhaustive} service, the visit period $V_i$ can be divided into the following subperiods: $V_i = V_i^{(S_i)}+V_i^{(V_{i+1})}+\dots+V_i^{(S_{i+N-1})}+V_i^{(V_{i})}$. First the type $i^{(S_i)}$ customers that were present at the visit beginning are served, followed by the type $i^{(V_{i+1})}$ customers, and so on. Note that during these services only type~$j^{(V_i)}$ customers arrive in $Q_j$, $j=1,\dots,N$. Visit period $V_i$ ends with $V_i^{(V_i)}$, i.e., the exhaustive service of all type $i^{(V_i)}$ customers that have arrived during $V_i$ so far.
The joint queue length process at polling instants of each of the subperiods is still a Multi-Type Branching Process, because the service discipline in each subperiod satisfies the Branching Property.
Hence, the laws of motion can be obtained by applying this property successively. As an example, we show the relations for the PGFs of the joint queue length distributions at beginnings and endings of the subperiods of $V_1$:
\begin{align*}
\VB_{1}^{(V_2)}(\z) = \VC_{1}^{(S_1)}(\z)&=\VB_{1}^{(S_1)}\left(z_1^{(V_1)},\beta_1\big(\sum_{j=1}^N\lambda_j^{(V_1)}(1-z_j^{(V_1)})\big),z_1^{(V_2)},\dots,z_N^{(S_N)}\right),\\
\VB_{1}^{(S_2)}(\z)=\VC_{1}^{(V_2)}(\z)&=\VB_{1}^{(V_2)}\left(z_1^{(V_1)},1,\beta_1\big(\sum_{j=1}^N\lambda_j^{(V_1)}(1-z_j^{(V_1)})\big),z_1^{(S_2)},\dots,z_N^{(S_N)}\right),\\
\vdots\quad &\\
\VB_{1}^{(V_1)}(\z)=\VC_{1}^{(S_N)}(\z)&=\VB_{1}^{(S_N)}\left(z_1^{(V_1)},1,\dots,1,\beta_1\big(\sum_{j=1}^N\lambda_j^{(V_1)}(1-z_j^{(V_1)})\big),z_2^{(V_1)},\dots,z_N^{(S_N)}\right),\\
\VC_{1}^{(V_1)}(\z)&=\VB_{1}^{(V_1)}\left(\pi_1\big(\sum_{j\neq1}\lambda_j^{(V_1)}(1-z_j^{(V_1)})\big),1,\dots,1,z_2^{(V_1)},\dots,z_N^{(S_N)}\right).
\end{align*}
During a switch-over time $S_j$ only type $i^{(S_j)}$ customers arrive, $i,j=1,\dots,N$. We can relate the PGF of the joint queue length distribution at the beginning of a visit to $Q_2$ (starting with the service of type $2^{(S_2)}$ customers) to $\VC_{1}^{(V_1)}(\z)$:
\[
\VB_{2}^{(S_2)}(\z)=\VC_{1}^{(V_1)}(\z)\,\sigma_1\Big(\sum_{j=1}^N\lambda_j^{(S_1)}(1-z_j^{(S_1)})\Big).
\]
%Iteration of this expression leads to the steady-state joint queue length distribution at beginnings and completions of all subperiods within each visit period.
The above expressions can be used to express $\VB_{2}^{(S_2)}(\cdot)$ in terms of $\VB_{1}^{(S_1)}(\cdot)$, and this can be repeated to obtain a recursion for $\VB_{1}^{(S_1)}(\cdot)$.

\begin{remark}
For \emph{gated} service we take similar steps, but they are slightly different because arriving customers will always be served in the next cycle. This means that a visit to $Q_i$ starts with the service of all type $i^{(V_i)}$ customers present at that polling instant: $V_i = V_i^{(V_{i})} + V_i^{(S_i)}+V_i^{(V_{i+1})}+\dots+V_i^{(S_{i+N-1})}$. The relations for the PGF of the joint queue length distribution at beginnings and endings of the subperiods of $V_1$ are:
\begin{align*}
\VB_{1}^{(S_1)}(\z) = \VC_{1}^{(V_1)}(\z)&=\VB_{1}^{(V_1)}\left(\beta_1\big(\sum_{j=1}^N\lambda_j^{(V_1)}(1-z_j^{(V_1)})\big),z_1^{(S_1)},\dots,z_N^{(S_N)}\right),\\
\VB_{1}^{(V_2)}(\z) = \VC_{1}^{(S_1)}(\z)&=\VB_{1}^{(S_1)}\left(z_1^{(V_1)},\beta_1\big(\sum_{j=1}^N\lambda_j^{(V_1)}(1-z_j^{(V_1)})\big),z_1^{(V_2)},\dots,z_N^{(S_N)}\right),\\
\vdots\quad &\\
\VC_{1}^{(S_N)}(\z)&=\VB_{1}^{(S_N)}\left(z_1^{(V_1)}, 1, \dots, 1,\beta_1\big(\sum_{j=1}^N\lambda_j^{(V_1)}(1-z_j^{(V_1)})\big),z_2^{(V_1)},\dots,z_N^{(S_N)}\right).
\end{align*}
\end{remark}
The remainder of this section is valid for any branching-type service discipline treating customers in order of arrival in each queue, such as, e.g., exhaustive, gated, globally gated and multi-stage gated \cite{vdmeiresing07}.
Having determined the joint queue length distribution at beginnings and endings of all subperiods within each visit period, we are ready to determine the joint queue length distribution at departure epochs of all customer subtypes. We follow the approach in~\cite{semphd,borst97}, which itself is based on Eisenberg's approach \cite{eisenberg72}, developing a relation between joint queue lengths at service  beginnings/completions and visit beginnings/endings. In \cite{semphd}, for conventional polling systems, the joint distribution of queue length vector and server position at service completions leads to the marginal queue length distribution. Developing an equivalent for our model, requires distinguishing between customer subtypes. Firstly, the queue length vector $\z$ contains all customer subtypes. Secondly, the type of service completion is not just defined by the location $i$ of the server, but also by the subtype $P$ of the customer that has been served. Therefore, let $M_i^{(P)}(\z)$ denote the PGF of the joint distribution of the subtypes of customers being served (combination of $i=1,\dots,N$ and $P \in \{V_1, S_1, \dots, V_N, S_N\}$) and queue length vector of all customer subtypes at service completions.
Equation (3.4) in \cite{semphd}, applied to our model, gives:
\begin{equation}
M_i^{(P)}(\z) = \frac{1}{\overline{\lambda}\E[C]}\,\frac{\beta_i\left(\sum_{j=1}^N\lambda_j^{(V_i)}(1-z_j^{(V_i)})\right)}{z_i^{(P)}-\beta_i\left(\sum_{j=1}^N\lambda_j^{(V_i)}(1-z_j^{(V_i)})\right)}
\left[ \VB_{i}^{(P)}(\z) - \VC_{i}^{(P)}(\z)\right], \label{jointQLandServerPositionAtdeparture}
\end{equation}
for $i=1,\dots,N; P \in \{V_1, S_1, \dots, V_N, S_N\}$, and $\overline{\lambda} = \sum_{i=1}^N \overline{\lambda}_i$. %We now focus on departures from an arbitrary queue, say $Q_i$, only.
Thus, $M_i^{(P)}(\z)$ is the generating function of the probabilities that, at an arbitrary departure epoch, the departing customer is a type $i^{(P)}$ customer \emph{and} the number of customers left behind by this departing customer is $l_1^{(V_1)}, \dots, l_N^{(S_N)}$.
We now focus on the queue length vector of subtypes of type $i$ customers only, \emph{given that} the departure takes place at $Q_i$. The probability that an arbitrary service completion (regardless of the subtype of the customer) takes place at $Q_i$, is $\overline{\lambda}_i/\overline{\lambda}$. It is convenient to introduce the notation $\zi = (1, \dots, 1, z_i^{(V_1)}, \dots, z_i^{(S_N)}, 1, \dots, 1)$. The PGF of the joint queue length distribution of all subtypes of type $i$ customers at an arbitrary departure from $Q_i$ is:
\begin{equation}
\E\left[\left(z_i^{(V_1)}\right)^{D_i^{(V_1)}}\cdots \left(z_i^{(S_N)}\right)^{D_i^{(S_N)}} \right] = \frac{\overline{\lambda}}{\overline{\lambda}_i}\sum_{j=1}^N \left(M_i^{(V_j)}(\zi)+M_i^{(S_j)}(\zi)\right)
\label{jointQLatdeparture}
\end{equation}
where $D_i^{(P)}$ is the number of type $i^{(P)}$ customers left behind at a departure from $Q_i$ (which should not be confused with $L_i^{(P)}$, the number of type $i$ customers at an arbitrary moment while the server is at position $P$).

\begin{remark}
Substitution of $z_i^{(P)} = z$ for all $P\in \{V_1, S_1, \dots, V_N, S_N\}$ in \eqref{jointQLatdeparture} gives the marginal queue length distribution of type $i$ customers at departure epochs, which is equal to \eqref{queuelengthGFarrival}, the marginal queue length distribution at arrival epochs of a type $i$ customer.
\end{remark}

\subsection{Waiting times}

Now we present a generalisation of the distributional form of Little's law that can be applied to the joint queue length distribution of all subtypes of a type $i$ customer at departure epochs from $Q_i$, to obtain the waiting time LST of a type $i$ customer.

\begin{theorem}\label{distformLittleTheorem}
The LST of the distribution of the waiting time $W_i$ of a type $i$ customer, $i=1,\dots,N$, is given by:
\begin{equation}
\E\left[\ee^{-\omega W_i}\right] = \frac{1}{\beta_i(\omega)}\E\left[\left(1-\frac{\omega}{\lambda_i^{(V_1)}}\right)^{D_i^{(V_1)}}\cdots \left(1-\frac{\omega}{\lambda_i^{(S_N)}}\right)^{D_i^{(S_N)}} \right].
\label{waitingtimeLST}
\end{equation}
\begin{proof}
We focus on the departure of a %particular
type $i$ customer that arrived during $P_A\in \{V_1, S_1, \dots, V_N, S_N\}$. We make use of the fact that the sojourn time (i.e., waiting time plus service time) of this tagged type~$i^{(P_A)}$ customer can be determined by studying the subtypes of all type $i$ customers that he leaves behind on his departure.
We need to distinguish between two cases, which can be treated simultaneously, but require different notations.
Firstly, the case where a customer arrives in the system and departs during another period. In the second case, the customer departs during the same period in which he arrived. Obviously, in our model the second case can only occur if a customer arrives at a queue with exhaustive service while it is being visited by the server.

\paragraph{Case 1: departure in a different period.}
In this case we have that $P_A \neq V_i$, or $P_A = V_i$ but the cycle in which the arrival took place is not the same as the cycle in which the departure takes place %(a situation that might occur, e.g., with gated service).
(this situation cannot occur with exhaustive service).
All type $i$ customers that are left behind, have arrived during the residual period $P_A$, all periods between $P_A$ and $V_i$ (if any), and during the elapsed part of $V_i$. Denote by $P_I$ the set of visit periods and switch-over periods that lie between $P_A$ and $V_i$. Furthermore, let $P_{A,\textrm{res}}$ be the residual period $P_A$. % of the type $i^{(P_A)}$ customer. %Finally denote by $V_{i,\textrm{past}}$ the elapsed time of $V_i$ at the departure instant of a randomly tagged type $i$ customer.
Finally denote by $V_{i,\textrm{past}}$ the age of $V_i$ at the departure instant of the tagged type $i$ customer.

\paragraph{Case 2: departure during the period of arrival.}
If the customer arrived during the same visit period in which his departure takes place, %i.e. $P_A = V_i$,
take $P_{A,\textrm{res}}=0, P_I=\emptyset$, and $V_{i,\textrm{past}}$ is the time that elapsed since the arrival of the tagged type $i^{(V_i)}$ customer.

In both cases, the joint queue length distribution of all customer $i$ subtypes at this departure instant is given by \eqref{jointQLatdeparture}.
%Note that at such a departure instant, there are no type $i$ customers present anymore that have arrived during the visit periods and switch-over periods between the previous $V_i$ and arrival period $P_A$. This results in:
Since we assume FCFS service, at such a departure instant no type $i$ customers are present anymore that have arrived before the arrival epoch of the tagged type $i$ customer. This results in:
\begin{equation}
\E\left[\left(z_i^{(V_1)}\right)^{D_i^{(V_1)}}\cdots \left(z_i^{(S_N)}\right)^{D_i^{(S_N)}} \right] = \E\left[\ee^{-\lambda_i^{(P_A)}(1-z_i^{(P_A)})P_{A,\textrm{res}}- \sum_{p\in P_I}\lambda_i^{(p)}(1-z_i^{(p)})p-\lambda_i^{(V_i)}(1-z_i^{(V_i)})V_{i,\textrm{past}}}\right].
\label{relationLandW}
\end{equation}
Equation \eqref{waitingtimeLST} follows from the relation $W_i + B_i = P_{A,\textrm{res}} + \sum_{p\in P_I}p + V_{i,\textrm{past}}$ and substitution of $z_i^{(P)} = 1-\frac{\omega}{\lambda_i^{(P)}}$ for all $P\in \{V_1, S_1, \dots, V_N, S_N\}$ in \eqref{relationLandW}.
\end{proof}
\end{theorem}

\begin{remark}
Theorem \ref{distformLittleTheorem} only holds if $\lambda_i^{(P)}>0$ for all $i=1,\dots,N$, and $P \in \{V_1, S_1, \dots, V_N, S_N\}$. If $\lambda_i^{(P)} = 0$ for a certain $i$ and $P$, we can still find an expression for $\E\left[\ee^{-\omega W_i}\right]$, but we might have to resort to some ``tricks''. In Section \ref{examples}, Example 2, we show how the introduction of an extra (virtual) customer type can help to resolve this problem.

%Since \eqref{jointQLandServerPositionAtdeparture} gives the joint distribution of all customer subtypes, it is possible to use $\E\left[\left(z_i^{(V_1)}\right)^{N_i^{(V_1)}}\cdots \left(z_j^{(P)}\right)^{N_j^{(P)}}\cdots\left(z_i^{(S_N)}\right)^{N_i^{(S_N)}} \right]$ instead of \eqref{jointQLatdeparture} and adapt \eqref{waitingtimeLST} accordingly. The only restriction that applies to which customer types $j^{(P)}$ are suitable, is the requirement that this type should not be served before the type $i$ customer, because we need these customers to be in the system at the moment of departure of the specified type $i^{(P)}$ customer.
\end{remark}

%\begin{remark}
%Theorem \ref{distformLittleTheorem} is valid for any service discipline that only allows a customer to receive service, after all customers in the same queue that have arrived earlier, have been served. In particular, it is valid for gated and exhaustive service, even though for gated service a departing type $i^{(V_i)}$ customer might leave behind type $i^{(V_i)}$ customers that have arrived during the cycle in which he arrived, plus type $i^{(V_i)}$ customers that have arrived during the cycle in which he departs. For non-branching service disciplines, one generally cannot determine the joint queue length distribution at departure epochs \eqref{jointQLatdeparture}, which makes it impossible to compute \eqref{waitingtimeLST} explicitly.
%\end{remark}

\section{Cycle time, intervisit time and visit times\label{cycletimesection}}

In the previous sections we repeatedly needed the mean cycle time $\E[C]$ and the mean visit times $\E[V_i],\ i=1,\dots,N$. In this section we study the LSTs of the cycle time distribution and visit time distributions, which can be used to obtain the mean and higher moments. The LSTs of the distributions of the visit times $V_i,\ i=1,\dots,N$, can easily be determined for any branching-type service discipline using the function $\theta_i(\cdot)$, introduced in Remark \ref{thetadef}, and the joint queue length distribution at the visit beginning of $Q_i$ (not taking subtypes into account):
\begin{equation}
\E[\ee^{-\omega V_i}] = \LB^{(V_{i})}(1,\dots,1,\theta_i(\omega),1, \dots, 1).\label{LSTVi}
\end{equation}
The cycle time $C_i$ is defined as the time that elapses between two consecutive visit beginnings to $Q_i$. We consider branching-type service disciplines only, i.e., service disciplines for which Property \ref{resingproperty} holds. The cycle time LST for polling models with branching-type service disciplines and arrival rates independent of the server position, has been established in \cite{boxmafralixbruin08}. We adapt their approach to the model with arrival rates depending on the server location. Using $\theta_i(\cdot)$, $i=1,\dots,N$, we define the following functions in a recursive way:
\begin{align*}
\psi^{(V_N)}(\omega) &= \omega,\\
\psi^{(V_i)}(\omega) &= \omega + \sum_{k=i+1}^N\lambda_k^{(V_i)}\left(1-\theta_k(\psi^{(V_k)}(\omega))\right), \qquad i=N-1,\dots,1.
\end{align*}
Similarly, define:
\begin{align*}
\psi^{(S_N)}(\omega) &= \omega,\\
\psi^{(S_i)}(\omega) &= \omega + \sum_{k=i+1}^N\lambda_k^{(S_i)}\left(1-\theta_k(\psi^{(V_k)}(\omega))\right), \qquad i=N-1,\dots,1.
\end{align*}

\begin{theorem}
The LST of the distribution of the cycle time $C_1$ is:
\begin{equation}
\E\left[\ee^{-\omega C_1}\right] = \LB^{(V_{1})}\left(\theta_1(\psi^{(V_1)}(\omega)), \dots, \theta_N(\psi^{(V_N)}(\omega))\right)\,\prod_{i=1}^N \sigma_i\left(\psi^{(S_i)}(\omega)\right).
\label{cycletimeLST}
\end{equation}
\begin{proof}
Similar to the proof of Theorem 3.1 in \cite{boxmafralixbruin08}, by giving an expression for the cycle time LST conditioned on the numbers of customers in all queues at the beginning of a cycle, and then by subsequently unconditioning one queue at a time.
\end{proof}
\end{theorem}
The LST of the distribution of the intervisit time $I_1$ can be found in a similar way:
\begin{equation}
\E\left[\ee^{-\omega I_1}\right] = \LB^{(S_{1})}\left(1,\theta_2(\psi^{(V_2)}(\omega)), \dots, \theta_N(\psi^{(V_N)}(\omega))\right)\,\prod_{i=1}^N \sigma_i\left(\psi^{(S_i)}(\omega)\right).
\label{intervisittimeLST}
\end{equation}
Equations \eqref{cycletimeLST} and \eqref{intervisittimeLST} hold for general branching-type service disciplines. For gated and exhaustive service we can give expressions that are more compact and easier to interpret, using the joint queue length distribution of all customer subtypes at visit beginnings, as given in Subsection \ref{jointsubtypes}.
\begin{theorem}
If $Q_i$ receives \emph{exhaustive service}, the LST of the distribution of the cycle time $C_i^*$, starting at a visit \emph{ending} to $Q_i$, and the LST of the distribution of the intervisit time $I_i$, are given by:
\begin{align}
\E\big[\ee^{-\omega C^*_i}\big] &= \VB_{i}^{(S_i)}\big(1, \dots, 1, \pi_i(\omega)-\frac{\omega}{\lambda_i^{(V_1)}}, \dots, \pi_i(\omega)-\frac{\omega}{\lambda_i^{(S_N)}}, 1, \dots, 1\big),\label{cycletimeLSTexh}\\
\E\big[\ee^{-\omega I_i}\big] &= \VB_{i}^{(S_i)}\big(1, \dots, 1, 1-\frac{\omega}{\lambda_i^{(V_1)}}, \dots, 1-\frac{\omega}{\lambda_i^{(S_N)}}, 1, \dots, 1\big),\label{intervisittimeLSTexh}
\end{align}
provided that $\lambda_i^{(P)} \neq 0$ for all $P \in \{V_1, S_1, \dots, V_N, S_N\}$. In the right hand sides of \eqref{cycletimeLSTexh} and \eqref{intervisittimeLSTexh}, the components $z_j^{(P)}$ with $j \neq i$ are 1.

If $Q_i$ receives \emph{gated service}, the LST of the distribution of the cycle time $C_i$, and the LST of the distribution of the intervisit time $I_i$, are given by:
\begin{align}
\E\big[\ee^{-\omega C_i}\big] &= \VB_{i}^{(V_i)}\big(1, \dots, 1, 1-\frac{\omega}{\lambda_i^{(V_1)}}, \dots, 1-\frac{\omega}{\lambda_i^{(S_N)}}, 1, \dots, 1\big),\nonumber\\
%\E\big[\ee^{-\omega I_i}\big] &= \VB_{i}^{(V_i)}\big(1, \dots, 1, 1, 1-\frac{\omega}{\lambda_i^{(S_1)}}, \dots, 1-\frac{\omega}{\lambda_i^{(S_N)}}, 1, \dots, 1\big),\label{lstIgated}
\E\big[\ee^{-\omega I_i}\big] &= \VB_{i}^{(V_i)}\big(1, \dots, 1, 1, 1-\frac{\omega}{\lambda_i^{(V_1)}}, \dots, 1-\frac{\omega}{\lambda_i^{(S_{i-1})}},1,1-\frac{\omega}{\lambda_i^{(S_i)}},\dots,1-\frac{\omega}{\lambda_i^{(S_N)}}, 1, \dots, 1\big),\label{lstIgated}
\end{align}
again provided that $\lambda_i^{(P)} \neq 0$ for all $P \in \{V_1, S_1, \dots, V_N, S_N\}$. Note that $z_i^{(V_i)} = 1$ in \eqref{lstIgated}.
\begin{proof}
We prove the exhaustive case only, the proof for gated service proceeds along the same lines. Using $I_i = S_i + V_{i+1}+ S_{i+1}+\dots+S_{i+N-1}$, and the fact that no type $i^{(V_i)}$ customers are present at the beginning of the intervisit period (and hence also at the beginning of a cycle $C_i^*$), we obtain:
\begin{equation}
\VB_{i}^{(S_i)}\left(1, \dots, 1, z_i^{(V_1)}, \dots, z_i^{(S_N)}, 1, \dots, 1\right) = \E\left[\ee^{-\lambda_i^{(S_i)}(1-z_i^{(S_i)})S_i- \dots-\lambda_i^{(S_{i+N-1})}(1-z_i^{(S_{i+N-1})})S_{i+N-1}}\right].\label{relationvbandcycletime}
\end{equation}
Substitution of $z_i^{(P)} = 1-\frac{\omega}{\lambda_i^{(P)}}$ for all $P \in \{V_1, S_1, \dots, V_N, S_N\}$ proves \eqref{intervisittimeLSTexh}. Equation \eqref{cycletimeLSTexh} follows by using the relation $C_i^* = I_i +V_i$, and noting that $V_i$ is the sum of the busy periods initiated by all type $i$ customers that have arrived during $I_i$. In terms of LSTs:
\begin{align*}
\E\left[\ee^{-\omega C_i^*}\right] &= \E\left[\ee^{-\left(\omega+\lambda_i^{(S_i)}(1-\pi_i(\omega))\right)S_i- \dots-\left(\omega+\lambda_i^{(S_{i+N-1})}(1-\pi_i(\omega))\right)S_{i+N-1}}\right]\\
&=\E\left[\ee^{-\lambda_i^{(S_i)}\left(1-\big(\pi_i(\omega)-\frac{\omega}{\lambda_i^{(S_i)}}\big)\right)S_i- \dots-\lambda_i^{(S_{i+N-1})}\left(1-\big(\pi_i(\omega)-\frac{\omega}{\lambda_i^{(S_{i+N-1})}}\big)\right)S_{i+N-1}}\right],
\end{align*}
which, by \eqref{relationvbandcycletime}, reduces to \eqref{cycletimeLSTexh}.
\end{proof}
\end{theorem}
Differentiation of the LSTs of $C_i$ and $C_i^*$ for $i=1,\dots,N$, shows that, just like in polling models with constant arrival rates, the mean cycle time does not depend on the starting point of the cycle, i.e. $\E[C_i]=\E[C_i^*]=\E[C]$. The mean cycle time $\E[C]$ and mean visit times $\E[V_i]$ can be obtained by differentiating the corresponding LSTs. In the next section a more efficient method is described to compute them.

\section{Mean Value Analysis\label{mva}}

In this section we extend the Mean Value Analysis (MVA) framework for polling models, originally developed by Winands et al.~\cite{winands06}, to suit the concept of smart customers. For this purpose, we first outline the main ideas of MVA for polling systems. Subsequently, we determine the mean visit times and the mean cycle time in a numerically more efficient way than in the previous section, and, finally, we present the MVA equations for a polling system with smart customers.

\subsection{Main idea MVA}

For ``ordinary'' polling models, where the arrival rates at a queue do not depend on the position of the server, in~\cite{winands06} an approach is described for deriving the steady-state mean waiting times at each of the queues, $\E[W_i]$ for $i=1,\dots,N$, by setting up a system of linear equations, where each equation has a probabilistic and intuitive explanation. We sketch the main ideas of MVA for exhaustive service; the cases of gated or mixed service disciplines require only minor changes.

The mean waiting time $\E[W_i]$ of a type~$i$ customer, excluding his service time, can be expressed in the following way: upon arrival of a (tagged) type~$i$ customer, he has to wait for the (remaining) time it takes to serve all type~$i$ customers already present in the system, plus possibly the time before the server arrives at~$Q_i$. By PASTA, the arriving customer finds in expectation $\E[\LQ_i]$ waiting type~$i$ customers in queue, each having an expected service time $\E[B_i]$. Note that we use $\LQ_i$ to denote the queue lengths \emph{excluding} customers in service. The expected time until the server returns to~$Q_i$, is denoted by $\E[T_i]$ (which depends on the service discipline of all queues).
A fraction $\rho_i := \lambda_i \, \E[B_i]$ of the time, the server is serving~$Q_i$, and hence, with probability $\rho_i$, an arriving customer has to wait for a mean residual service time, denoted by $\E[R_{B_i}]$; otherwise he has to wait until the server returns. This gives, for $i=1,\dots,N$:
\begin{equation*}\label{eq:EwiOrdinaryMVA}
\E[W_i] = \E[\LQ_i]\, \E[B_i] + \rho_i\, \E[R_{B_i}] + (1 - \rho_i)\, \E[T_i].
\end{equation*}
Little's law gives $\E[\LQ_i] = \lambda_i\, \E[W_i]$, for $i=1,\dots,N$, and so it remains to derive $\E[T_i]$. For this, first a system of equations is composed for the \emph{conditional} mean queue lengths, which can be expressed in mean residual durations of (sums of) visit and switch-over times. The solution of this system of equations can be used to determine $\E[T_i]$, and hence $\E[\LQ_i]$ and $\E[W_i]$ follow.

\subsection{Mean visit times and mean cycle time}

For the case of smart customers, the visit times to a queue depend on all arrival rates $\lambda_i^{(V_j)}$ and $\lambda_i^{(S_j)}$. In order to extend MVA to this case, we first derive the mean visit times to each of the queues, $\E[V_i]$, for $i=1,\ldots,N$. We set up a system of~$N$ linear equations where the mean visit time of a queue is expressed in terms of the other mean visit times. We again focus on the exhaustive service discipline.

At the moment the server finishes serving $Q_i$, there are no type~$i$ customers present in the system any more. From this point on, the number of type~$i$ customers builds up at rates $\lambda^{(S_{i})}, \lambda^{(V_{i+1})}, \dots, \lambda^{(S_{i+N-1})}$ (depending on the position of the server), until the server starts working on~$Q_i$ again. Each of these customers initiates a busy period, with mean $\E[\BP_i] := \E[B_i]/(1-\lambda_i^{(V_i)} \E[B_i])$. This gives:
\begin{equation*}
\E[V_i] = \E[\BP_i] \left(\lambda_i^{(S_i)} \E(S_i) + \sum_{j=i+1}^{i+N-1} \left(\lambda_i^{(V_j)} \E[V_j] + \lambda_i^{(S_j)} \E[S_j] \right)\right),
%\label{eqmva_EVi}
\end{equation*}
for $i=1,\ldots,N$. % where all summations should be understood as to be cyclically, and
The~$\E[V_i]$ follow from solving this set of equations. This method is computationally faster than determining (and differentiating) the LSTs of the visit time distributions \eqref{LSTVi}. Once the mean visit times have been obtained, the mean cycle time follows from $\E[C] = \sum_{i=1}^N(\E[V_i]+\E[S_i])$.

\subsection{MVA equations}

We extend the MVA approach to polling systems with smart customers. First, we briefly introduce some extra notation, then we give expressions for the mean waiting times, and the mean conditional and unconditional queue lengths. %, and the residual interval lengths.
After eliminating variables, we end up with a system of %$(2N)^2$
linear equations. The system can (numerically) be solved in order to find the unknowns, in particular, the mean unconditional queue lengths and the mean waiting times. Although all equations are discussed in the present section, for the sake of brevity of this section, some of them are presented in Appendix \ref{appendixMVAeqns}.

%As the arrival rates differ during a visit time and switch-over time, we split a cycle as follows: $V_1, S_1, V_{2},$ $\ldots,$ $V_N, S_N$, where $V_i$ ($S_i$) denotes the visit time to (switch-over time from) $Q_i$, with mean duration $\E[V_i]$ ($\E[S_i]$). To both we refer as a \emph{period}, i.e.\ either $V_i$ or $S_i$.
The fraction of time the system is in a given period $P \in \{V_1, S_1, \dots, V_N, S_N\}$ is denoted by $q^{(P)}:= \frac{\E[P]}{\E[C]}$. The mean residual duration of a period $P$, at an arbitrarily chosen point in this period, is denoted by $\E[R_{P}] = \frac{\E[P^2]}{2\E[P]}$.
The mean conditional number of type~$j$ customers in the queue during a random point in $P$ is denoted by $\E[\LQ_j^{(P)}]$, and the mean (unconditional) number of type~$j$ customers in queue is denoted by $\E[\LQ_j]$. Note that $\LQ_j$ and $\LQ_j^{(P)}$ do \emph{not} include a potential customer in service, whereas $L_j$ and $L_j^{(P)}$, introduced in Section \ref{nosubtypes}, denote queue lengths including customers being served. %Throughout this section queue lengths include customers waiting in the queue only, excluding customers in service.

We define an \emph{interval}, e.g.\ $(V_i\!:\!S_j)$, as the consecutive periods from the first mentioned period on, until and including the last mentioned period, here consisting of the periods $V_i, S_i, V_{i+1}, S_{i+1},$ $\dots,$ $V_j, S_j$.
The mean residual duration of an interval, e.g.\ $(V_i\!:\!S_j)$, is denoted by $\E[R_{V_i : S_j}]$. Analogously, we define $\E[R_{V_i : V_j}]$, $\E[R_{S_i : V_j}]$ and $\E[R_{S_i : S_j}]$.

An important concept in the remainder of the analysis is the concept of \emph{conditional durations} of a period. This is an extension of the well-known residual duration, or the age of a period. It deals with the length of a period within the cycle (i.e., a visit time or a switch-over time), given that the system is being observed from another period. Before we proceed, we clarify this important concept by a simple example. Consider a vacation system, i.e., a polling system with $N=1$. A cycle consists of a switch-over time (or: vacation) $S_1$, followed by a visit time $V_1$, during which the queue is served exhaustively. Now assume that the system is being observed at a random epoch during the switch-over time $S_1$. We derive an expression for $\EintForw{V_1}{S_1}$, which is the conditional mean visit time \emph{following} the switch-over time, \emph{given that} the system is being observed during $S_1$. Since service is exhaustive, the visit time consists of the busy periods of the customers that arrived during the elapsed part of $S_1$, denoted by $S_{1,\textrm{past}}$, plus the busy periods of the customers that will arrive during the residual switch-over time, denoted by $S_{1,\textrm{res}}$. Hence, it can be seen that in this system
\[
\EintForw{V_1}{S_1} = \frac{\lambda_1^{(S_1)}\E[B_1]}{1-\lambda_1^{(V_1)}\E[B_1]}\left(\E[S_{1,\textrm{past}}]+\E[S_{1,\textrm{res}}]\right)
=\frac{\lambda_1^{(S_1)}\E[B_1]}{1-\lambda_1^{(V_1)}\E[B_1]}\frac{\E[S_1^2]}{\E[S_1]}.
\]
Instead of studying the mean visit time \emph{following} the switch-over time during which the system is observed, we can also study the mean visit time \emph{preceding} this particular switch-over time, denoted by $\EintBack{V_1}{S_1}$. Now the expression is easier, because a switch-over time is independent of the preceding visit time, so
\[
\EintBack{V_1}{S_1} = \E[V_1] = \frac{\lambda_1^{(S_1)}\E[B_1]}{1-\lambda_1^{(V_1)}\E[B_1]}\E[S_1].
\]
Similarly, we denote by $\EintForw{S_1}{V_1}$ and $\EintBack{S_1}{V_1}$ the mean switch-over times \emph{following}, respectively \emph{preceding}, the visit period $V_1$ during which the system is observed. Because of the independence between $S_1$ and the preceding $V_1$, it is immediately clear that
\[
\EintForw{S_1}{V_1} = \E[S_1].
\]
Unfortunately, $\EintBack{S_1}{V_1}$ cannot be determined this easily, since $V_1$ is positively correlated with the preceding switch-over time. However, it plays a role in the set of MVA equations developed later in this section, relating the conditional mean queue lengths, the conditional mean waiting times and the conditional durations of the periods in a cycle. We leave it up to the reader to verify that Equation \eqref{lemmaeqn} reduces to $\EintBack{S_1}{V_1} = \E[S_1^2]/\E[S_1]$ in the exhaustive vacation model.
This example simply serves the purpose of illustrating the concept of these conditional durations. In a polling system consisting of multiple queues, these expressions become more complicated and can only be found by solving sets of equations, as will be shown in the remainder of this section. Note that, because of conditional PASTA, an arbitrary customer arriving during $S_1$ finds the system in the same state as an observer who observes the system at an arbitrary epoch during $S_1$. Hence, the conditional durations of periods play an important role in determining the mean waiting times.

For the mean conditional durations of a period, we have the following: $\EintBack{V_i}{V_j}$ denotes the mean duration of the \emph{previous} period~$V_i$, seen from an arbitrary point in $V_j$ (i.e., $V_i$ seen backwards in time from the viewpoint of $V_j$), and $\EintForw{V_i}{V_j}$ denotes the mean duration of the \emph{next} period~$V_i$ (i.e., $V_i$ seen forwards in time from the viewpoint of $V_j$). For $i=j$ they both coincide, and represent the mean age, resp. the mean residual duration of $V_i$. Since the distribution of the age of a period is the same as the distribution of the residual period, we have $\EintBack{V_i}{V_i}=\EintForw{V_i}{V_i}=\E[R_{V_i}]$. Generally, however, $\EintBack{V_i}{V_j} \neq \EintForw{V_i}{V_j}$ for $i \neq j$, because of the dependencies between the durations of periods. Analogously, we define $\EintBack{V_i}{S_j}$, $\EintForw{V_i}{S_j}$, $\EintBack{S_i}{V_j}$ and $\EintForw{S_i}{V_j}$. Note that, e.g., $\EintForw{S_i}{V_j} = \E[S_i]$, but $\EintBack{S_i}{V_{j}} \neq \E[S_i]$.
As switch-over times are independent, the following quantities directly simplify:
\[
\EintBack{S_i}{S_j} = \EintForw{S_i}{S_j} =
\begin{cases}
\E[S_i] & \mbox{for $i \neq j$},\\
\E[R_{S_i}] & \mbox{for $i = j$}.
\end{cases}
\]
Having introduced the required notation, we now present the main theorem of this section, which gives a set of equations that can be solved to find the mean waiting times of customers in the system.

\begin{theorem}\label{mvatheorem}
The mean waiting times, $\E[W_i]$, for $i=1,\dots,N$, and the mean queue lengths, $\E[\LQ_i]$, satisfy the following equations:
\begin{align}
\E[W_i]&=\frac{q^{(V_i)} \lambda_i^{(V_i)}}{\overline{\lambda}_i}
         \left(\E[\LQ_i^{(V_i)}]  \E[B_i] +\E[R_{B_i}]\right)\nonumber\\
       & \quad \mbox{} + \sum_{j=i+1}^{i+N-1} \frac{q^{(V_j)} \lambda_i^{(V_j)}}{\overline{\lambda}_i}
         \left(\E[\LQ_i^{(V_j)}]  \E[B_i]+ \sum_{k=j}^{i+N-1}\left( \E[S_k] + \EintForw{V_k}{V_j} \right)\right)\nonumber\\
       & \quad \mbox{} + \sum_{j=i}^{i+N-1}
         \frac{q^{(S_j)} \lambda_i^{(S_j)}}{\overline{\lambda}_i}
         \left(\E[\LQ_i^{(S_j)}] \E[B_i] + \E[R_{S_j}] + \sum_{k=j+1}^{i+N-1} \left(\E[S_k] + \EintForw{V_k}{S_j} \right)\right),\label{eqmva_EWi1}\\
\E[\LQ_i] &= \overline{\lambda}_i \E[W_i],\label{eqmva_ELi_Little}\\
\E[\LQ_i] &= \sum_{j=i+1}^{i+N} \left( q^{(V_j)} \E[\LQ_i^{(V_j)}] + q^{(S_j)} \E[\LQ_i^{(S_j)}] \right),\label{eqmva_ELi_weightedsum}
\end{align}
where the conditional mean queue lengths $\E[\LQ_i^{(V_j)}]$ and $\E[\LQ_i^{(S_j)}]$, for $j=i+1,\dots,i+N-1$, are given by
\begin{align}
\E[\LQ_i^{(V_j)}] &= \sum_{k=i+1}^{j} \lambda_i^{(V_k)} \EintBack{V_k}{V_j} + \sum_{k=i}^{j-1} \lambda_i^{(S_k)} \EintBack{S_k}{V_j},
\label{eqmva_ELiVj}\\
\E[\LQ_i^{(S_j)}] &= \sum_{k=i+1}^{j} \lambda_i^{(V_k)}  \EintBack{V_k}{S_j} + \sum_{k=i}^{j} \lambda_i^{(S_k)} \EintBack{S_k}{S_j},\label{eqmva_ELiSj}
\end{align}
and where all $\EintBack{P_1}{P_2}$ and $\EintForw{P_1}{P_2}$, for $P_1, P_2 \in \{V_1, S_1, \dots, V_N, S_N\}$, satisfy the set of equations \eqref{eqmva_EViVif} -- \eqref{lemmaeqn} below, and \eqref{eqmva_ERSiSj_forw}--\eqref{eqmva_ERViSj} in Appendix \ref{appendixMVAeqns}.

\begin{proof}
In order to derive the mean waiting time $\E[W_i]$, we condition on the period in which a type~$i$ customer  arrives.
A fraction $q^{(V_j)} \lambda_i^{(V_j)}/ \overline{\lambda}_i $, and $q^{(S_j)} \lambda_i^{(S_j)}/ \overline{\lambda}_i $ respectively, of the type~$i$ customers arrives during period~$V_j$, and during period $S_j$ respectively. % (which are \emph{not} equal to $q^{(V_j)}$, $q^{(S_j)}$).
If a tagged type~$i$ customer arrives during period $V_i$ (i.e., while his queue is being served), he has to wait for a residual service time, plus the service times of all type~$i$ customers present in the system upon his arrival, which is (by conditional PASTA), $\E[\LQ_i^{(V_i)}]$. As a fraction $q^{(V_i)} \lambda_i^{(V_i)} / \overline{\lambda}_i$ of the customers arrives during $V_i$, this explains the first line of \eqref{eqmva_EWi1}. If the customer arrives in any other period, he has to wait until the server returns to $Q_i$ again. For this, we condition on the period in which he arrives. If the arrival period is a visit to~$Q_j$, say $V_j$ for $j \neq i$, he has to wait for the residual duration of $V_j$ and the interval $(S_j\!:\!S_{i-1})$, and for the service of the type~$i$ customers present in the system upon his arrival%, which is $\E[\LQ_i^{(V_j)}]$
. This gives the second line of \eqref{eqmva_EWi1}. The third line, the case where the customer arrives during the switch-over time from~$Q_j$ to $Q_{j+1}$ (period~$S_j$), can be interpreted along the same lines as the case~$V_j$.

Equation \eqref{eqmva_ELi_weightedsum} is obtained by unconditioning the conditional queue lengths $\E[\LQ_i^{(P)}]$.
%These conditional queue lengths can be expressed as follows.
The mean number of type~$i$ customers in the queue at an arbitrary point during~$V_j$, given by \eqref{eqmva_ELiVj}, is the mean number of customers built up from the last visit to~$Q_i$ (when~$Q_i$ became empty) until and including a residual duration of~$V_j$ (as the mean residual duration of $V_j$ is equal to the mean age of that period), taking into account the varying arrival rates. The mean number of type~$i$ customers queueing in the system during period~$S_j$, given by \eqref{eqmva_ELiSj}, can be found similarly.
Equations \eqref{eqmva_ELiVj} and \eqref{eqmva_ELiSj} show one of the difficulties in adapting the ``ordinary'' MVA approach to that of smart customers. If the arrival rates remain constant during a cycle, these expressions would reduce to $\lambda_i$ multiplied by the mean time passed since the server has left~$Q_i$. However, for the smart customers case, we have to keep track of the duration of all the intermediate periods, from the viewpoint of period~$V_j$ respectively~$S_j$.

As indicated in Theorem \ref{mvatheorem}, at this point, the number of equations is insufficient to find all the unknowns, $\EintBack{P_1}{P_2}$ and $\EintForw{P_1}{P_2}$, for $P_1, P_2 \in \{V_1, S_1, \dots, V_N, S_N\}$. In the remainder of the proof, we develop additional relations for these quantities to complete the set of equations.
We start by considering $\EintForw{V_i}{V_j}$, which is the mean duration of the next period~$V_i$, when observed from an arbitrary point in $V_j$. For $i=j$ this is just the residual duration of~$V_i$, consisting of a busy period induced by a customer with a residual service time left, and the busy periods of all type~$i$ customers in the queue. The cases $i \neq j$ need some more attention. The duration of~$V_i$ now consists of the busy period induced by the type~$i$ customers in the queue, which are in expectation $\E[\LQ_i^{(V_j)}]$ customers. During the periods $V_j, S_j, \ldots, S_{i-1}$, however, new type~$i$ customers are arriving, each contributing a busy period to the duration of~$V_i$. Hence, summing over these periods and taking into account the varying arrival rates, we get the mean total of newly arriving customers in this interval. This yields, for $i = 1,\dots,N$ and $j=i+1, \ldots, i+N-1$:
%eqns2c
\begin{align}
\EintForw{V_i}{V_i} & = \E[\BP_i] \,\E[\LQ_i^{(V_i)}] + \E[R_{B_i}]/\left(1-\lambda_i^{(V_i)} \E[B_i]\right),\label{eqmva_EViVif}
\\
\EintForw{V_i}{V_j} & = \E[\BP_i] \left(\E[\LQ_i^{(V_j)}]
                    + \sum_{k=j}^{i+N-1} \left(\lambda_i^{(V_k)} \EintForw{V_k}{V_j} +
                    \lambda_i^{(S_k)} \E[S_k]  \right)\right).\label{eqmva_EViVjf}
%                   \\& \hspace{5cm} \mbox{for $j=i+1, \ldots, i+N-1$.}\nonumber
\end{align}
Analogously $\EintForw{V_i}{S_j}$ denotes the mean duration of the next period~$V_i$, when observed from an arbitrary point in $S_j$. The explanation of its expression is along the same lines as that of $\EintForw{V_i}{V_j}$, although it should be noted that $i=j$ is not a special case. See \eqref{eqmva_EViSjf} in Appendix \ref{appendixMVAeqns}.

The last step in the proof of Theorem \ref{mvatheorem}, needs the following lemma to find the final relations between $\EintBack{P_1}{P_2}$ and $\EintForw{P_1}{P_2}$:
\begin{lemma}\label{lemmamva}
For $i = 1,\dots,N$, and $j = i+1,\dots,i+N$:
\begin{align}
\sum_{k=i}^{j-1} \frac{\E[S_k]}{\E[(S_i\!:\!V_j)]} &\left(\EintBack{S_i}{S_k} + \sum_{l=i+1}^k \left(\EintBack{S_l}{S_k} + \EintBack{V_l}{S_k} \right)\right.\nonumber\\
&\left.-\E[R_{S_k}] - \EintForw{V_j}{S_k} - \sum_{l=k+1}^{j-1}\left( \E[S_l] + \EintForw{V_l}{S_k}
\right)
\right)\nonumber\\
=
\sum_{k=i+1}^{j} \frac{\E[V_k]}{\E[(S_i\!:\!V_j)]}&\left( \EintForw{V_j}{V_k} +  \sum_{l=k}^{j-1}\left( \E[S_l] + \EintForw{V_l}{V_k} \right)\right.\nonumber\\
&\left.-\EintBack{S_i}{V_k} -  \EintBack{V_k}{V_k} - \sum_{l=i+1}^{k-1} \left(\EintBack{S_l}{V_k} + \EintBack{V_l}{V_k} \right)
\right).
\label{lemmaeqn}
\end{align}
\begin{proof}
Equation \eqref{lemmaeqn} can be proven by studying all mean residual interval lengths $\E[R_{S_i : V_j}]$, $\E[R_{S_i : S_j}]$, $\E[R_{V_i : V_j}]$ and $\E[R_{V_i : S_j}]$. %These can be expressed in two ways. For the first one, we consider the mean time the given interval still lasts. However, as the remaining duration of a random time interval equals in distribution the passed duration of it (and hence has equal mean), the other expression for the mean residual interval duration can be found by considering the time already passed since the beginning of the interval.
%
%Firstly, expressions are given for the duration of the time already passed during the interval, subsequently for the remaining duration.
Consider $\E[R_{S_i : V_j}]$, the mean residual duration of the interval $S_i, V_{i+1}, \ldots, V_j$. We condition on the period in which the interval is observed. As the mean duration of the interval is given by \mbox{$\E[(S_i\!:\!V_j)]$}, it follows that $\E[S_k]/\E[(S_i\!:\!V_j)]$ %(or $\E[V_k]/\E[(S_i\!:\!V_j)]$ respectively)
is the probability that the interval is observed in period~$S_k$% (or $V_k$ respectively)
. The remaining duration of the interval consists of the remaining duration of $S_k$ %(or $V_k$)
plus the mean durations of the (coming) periods $V_{k+1}, S_{k+1},\ldots, V_j$% (or $S_k, V_{k+1}, S_{k+1},\ldots, V_j$)
, when observed from period~$S_k$% (or $V_k$)
. When observing \mbox{$\E[(S_i\!:\!V_j)]$} from $V_k$, a similar way of reasoning is used.
This gives, for $i = 1,\dots,N$, and $j = i+1,\dots,i+N$:
%eqns33a
\begin{align}
\E[R_{S_i : V_j} ] &= \sum_{k=i}^{j-1} \frac{\E[S_k]}{\E[(S_i\!:\!V_j)]}\left( \E[R_{S_k}] + \EintForw{V_j}{S_k} + \sum_{l=k+1}^{j-1}\left( \E[S_l] + \EintForw{V_l}{S_k} \right)\right)\nonumber\\
                    &\quad \mbox{} + \sum_{k=i+1}^{j} \frac{\E[V_k]}{\E[(S_i\!:\!V_j)]}\left( \EintForw{V_j}{V_k} +  \sum_{l=k}^{j-1}\left( \E[S_l] + \EintForw{V_l}{V_k} \right)\right).
\label{eqmva_ERSiVj_forw}
\end{align}
We now use that the distribution of the residual length of an interval is the same as the distribution of the age of this interval.
Again, focus on $\E[R_{S_i : V_j}]$, conditioning on the period in which the interval is observed, but now looking forwards in time.
Consider all the periods in $(S_i : V_j)$ that have already passed when observing during $S_k$. %(or $V_k$)
The interval has lasted for the sum of these periods% (or $V_k$)
, plus the age of~$S_k$% (or $V_k$ respectively)
. The same can be done for an arbitrary point in $V_k$. This gives, for $i = 1,\dots,N$, $j = i+1,\dots,i+N$:
%
%eqns3a
\begin{align}
\E[R_{S_i : V_j} ] &= \sum_{k=i}^{j-1} \frac{\E[S_k]}{\E[(S_i\!:\!V_j)]}\left( \EintBack{S_i}{S_k} + \sum_{l=i+1}^k \left(\EintBack{S_l}{S_k} + \EintBack{V_l}{S_k} \right)\right)\nonumber\\
                    &\quad \mbox{} + \sum_{k=i+1}^{j} \frac{\E[V_k]}{\E[(S_i\!:\!V_j)]}\left( \EintBack{S_i}{V_k} +  \EintBack{V_k}{V_k} + \sum_{l=i+1}^{k-1} \left(\EintBack{S_l}{V_k} + \EintBack{V_l}{V_k} \right)\right).
\label{eqmva_ERSiVj}
\end{align}
The proof of Lemma \ref{lemmamva} is completed by equating \eqref{eqmva_ERSiVj_forw} and \eqref{eqmva_ERSiVj} and rearranging the terms.
\end{proof}
\end{lemma}
Similar to the proof of Lemma \ref{lemmamva}, we can develop two different expressions for each of the terms $\E[R_{S_i : S_j}], \E[R_{V_i : V_j}]$ and $\E[R_{V_i : S_j}]$. For the sake of brevity of this section, they are presented in Appendix \ref{appendixMVAeqns}, Equations \eqref{eqmva_ERSiSj_forw}--%\eqref{eqmva_ERViSj_forw} and \eqref{eqmva_ERSiSj}--
\eqref{eqmva_ERViSj}. Equating each pair of these expressions, completes the set of (linear) equations for the mean waiting times and mean queue lengths. This concludes the proof of Theorem \ref{mvatheorem}.
\end{proof}
\end{theorem}

\section{Pseudo-Conservation Law}\label{pclsection}

In this section we derive a so-called Pseudo-Conservation Law (PCL), which gives an expression for the weighted sum of the mean waiting times at each of the queues.
For ``ordinary'' cyclic polling systems, Boxma and Groenendijk~\cite{boxmagroenendijk87} derive a PCL under various service disciplines. This PCL, in commonly used notation $\rho_i=\lambda_i\E[B_i], \rho=\sum_{i=1}^N\rho_i, S=\sum_{i=1}^NS_i$, states that:
\begin{align}
\sum_{i=1}^N \rho_i \E[W_i]&=\rho \frac{\sum_{i=1}^N \rho_i \E[R_{B_i}]}{1-\rho}+\rho \E[R_S]
 +\frac{\E[S]}{2(1-\rho)} \left( \rho^2 - \sum_{i=1}^N \rho_i^2\right)+\sum_{i=1}^N\E[Z_{ii}],\label{eq:pclStandardExh}
\end{align}
with $Z_{ii}$ denoting the amount of work left behind by the server at $Q_i$ at the ending of a visit. For exhaustive service at $Q_i$, we have $\E[Z_{ii}]=0$, and for gated service $\E[Z_{ii}]=\frac{\rho_i^2\E[S]}{1-\rho}$.
%The PCL is derived by considering a workload decomposition, which we adapt to a model with smart customers.
%We focus on exhaustive service only. Other service disciplines, such as gated and $1$-limited, require marginal changes.

We base our approach on~\cite{boxmagroenendijk87}, and adapt their ideas to derive a PCL for a polling model with smart customers. The approach focusses on the mean amount of \emph{work} in the system at an arbitrary point in time. %, denoted by $\E[Y]$.
A required restriction for our approach in this section, is that the Poisson process according to which work arrives in the system, has a fixed arrival rate during all \emph{visit periods}. We also require that the amounts of work brought by an individual arrival are identically distributed for all visit periods.
We mention two typical cases where this requirement is satisfied. Firstly, the case when the arrival rate at a given queue stays constant during different \emph{visit} times, and secondly when the \emph{total} arrival rate remains constant during visit times \emph{and} the service times are identically distributed:
\begin{align}
\text{Case 1: }\qquad & \lambda_i^{(V_1)} = \lambda_i^{(V_2)} = \ldots = \lambda_i^{(V_N)}=:\lambda_i^{(V)},&&\qquad i=1,\dots,N,\label{eq:pclassump1}\\
\text{Case 2: }\qquad & \sum_{i=1}^N \lambda_i^{(V_j)} =: \Lambda^{(V)}, \text{ and }B_1 \equaldist \dots \equaldist B_N,  &&\qquad j=1,\dots,N.\label{eq:pclassump2}
\end{align}
Note that Case 1 does allow for different arrival rates during \emph{switch-over times}. During visit periods, let $\Lambda^{(V)}$ be the total arrival rate of all customer types, and let $B^{(V)}$ denote the generic service time of an arbitrary customer entering the system. In particular, this means for Case~1 that $\Lambda^{(V)} = \sum_{i=1}^N \lambda_i^{(V)}$ and $B^{(V)} \equaldist B_i$ with probability $\lambda_i^{(V)}/\Lambda^{(V)}$ for $i=1,\dots,N$.  We introduce $\rho^{(V)}$ to denote the mean amount of work entering the system per time unit during a visit period, so $\rho^{(V)} = \Lambda^{(V)} \E[B^{(V)}]$.

Denote by~$Y$ the amount of work in the polling system at an arbitrary point in time, and by~$Y^{(V)}$ and~$Y^{(S)}$
the amount of work at an arbitrary point during respectively a visit period, and a switch-over period.
Then
\begin{equation}
Y \equaldist \begin{cases}
	Y^{(V)}  & \mbox{w.p.\ } \overline{\rho},\\
	Y^{(S)}  & \mbox{w.p.\ } 1-\overline{\rho},
	\end{cases}
\label{eq:workdecomp}
\end{equation}
where $\overline{\rho} := \sum_{i=1}^N \overline{\rho}_i = \sum_{i=1}^N \overline{\lambda}_i \E[B_i]$ is the mean offered amount of work per time unit. Hence,
\begin{equation}
\E[Y] = \overline{\rho} \, \E[Y^{(V)}] + (1-\overline{\rho}) \E[Y^{(S)}].
\label{eq:EVdecomp}
\end{equation}
Another way to obtain the mean total amount of work in the system, is by taking the sum of the mean workloads. The mean workload in $Q_i$ is the mean amount of work of all customers in the queue, plus, with probability~$\overline{\rho}_i= \overline{\lambda}_i \E[B_i]$, the mean remaining amount of work of a customer in service at~$Q_i$: %
\begin{equation}
\E[Y] = \sum_{i=1}^N \left(\E[\LQ_i] \E[B_i] + \overline{\rho}_i \E[R_{B_i}] \right).\label{EYeqn2}
\end{equation}
In the next subsections we show that equating \eqref{eq:EVdecomp} and \eqref{EYeqn2}, and applying Little's law, $\E[\LQ_i] = \overline{\lambda}_i \E[W_i]$, gives a PCL for the mean waiting times in the system.
%
%In the sequel, we derive expressions for $\E[Y^{(V)}]$ and~$\E[Y^{(S)}]$.
But first we have to find $\E[Y^{(V)}]$ and $\E[Y^{(S)}]$. We start with the latter.

\subsection{Work during switch-over periods}

The term $\E[Y^{(S)}]$ denotes the mean amount of work in the system when observed at a random point in a switch-over interval. Denoting by $\E[Y^{(S_i)}]$ the mean amount of work in the system at an arbitrary moment during $S_i$, we can condition on the switch-over interval in which the system is observed:
\begin{equation}
\E[Y^{(S)}] = \sum_{i=1}^N \frac{\E[S_i]}{\E[S]} \E[Y^{(S_i)}].
\label{eq:EYinEYSj}
\end{equation}
We can split $\E[Y^{(S_i)}]$ into two parts: the mean amount of work present at the start of $S_i$, plus the mean amount of work built up since the start of the switch-over time. In expectation, a duration $\E[R_{S_i}]$ has passed since the beginning of the switch-over time, in which work arrived at rate $\lambda_j^{(S_i)} \E[B_j]$ at~$Q_j$. Hence, this gives a contribution to $\E[Y^{(S_i)}]$ of $\sum_{j=1}^N \lambda_j^{(S_i)} \E[B_j] \E[R_{S_i}]$. %Denote by $Z_{ii}$ the amount of work, present at $Q_i$ at the end of the visit to this queue.
For the work present at the start of the switch-over period, we start looking at the moment that the server left $Q_j$, leaving a mean amount of work $\E[Z_{jj}]$ behind in this queue.
For exhaustive service, $\E[Z_{jj}]=0$, for gated service $\E[Z_{jj}]= \lambda_j^{(V_j)}\E[B_j]\E[V_j]$.
Since then, the interval $(S_j\!:\!V_{i+N})$ has passed, for $j=i+1, \ldots, i+N-1$. % (as for $j=i$ the contribution is zero).
In this interval the amount of type~$j$ work increased at rates %$\lambda_j^{(S_j)} \E[B_j], \lambda_j^{(V)} \E[B_j], \ldots, \lambda_j^{(S_{i-1})} \E[B_j], \lambda_j^{(V)} \E[B_j]$
$\lambda_j^{(S_j)} \E[B_j], \lambda_j^{(V_{j+1})} \E[B_j], \ldots, \lambda_j^{(S_{i-1})} \E[B_j], \lambda_j^{(V_{i})} \E[B_j]$ during the various periods. This leads to the following expression for $\E[Y^{(S_i)}]$:
\begin{align}
\E[Y^{(S_i)}]
&= \sum_{j=1}^N \left(\lambda_j^{(S_i)} \E[B_j] \E[R_{S_i}]
              +\E[Z_{jj}]\right)
              +  \sum_{j=i+1}^{i+N-1} %\left(\E[Z_{jj}]+
                  \sum_{k=j}^{i+N-1} \left( \lambda_j^{(S_k)} \E[B_j]\E[S_k]
                            + \lambda_j^{(V_{k+1})} \E[B_j]\E[V_{k+1}] \right).\label{eq:EYSi}
\end{align}

\subsection{Work during visit periods}

The key observation in the proof of \cite{boxmagroenendijk87} is the \emph{work decomposition} property in a polling system. This property states that the amount of work at an arbitrary epoch in a visit period is distributed as the sum of two independent random variables: the amount of work in the ``corresponding'' $M/G/1$ queue at an arbitrary epoch during a busy period,  denoted by $Y_{M/G/1}^{(\textit{V})}$, and the amount of work in the polling system at an arbitrary epoch during a switch-over time, $Y^{(S)}$. In a polling model with smart customers, this decomposition does not typically hold, but a minor adaptation is required. We follow the proof in \cite{boxmagroenendijk87} as closely as possible, meaning that we use the concepts of ancestral line and offspring of a customer, as introduced in \cite{fuhrmanncooper85}.  We also copy the idea of comparing the polling system to an $M/G/1$ queue with vacations and Last-Come-First-Served (LCFS) service. The traffic process offered to this $M/G/1$ queue is identical to the traffic process of the polling system. The server of the $M/G/1$ queue takes vacations exactly during the switching periods of the polling system. These vacations might interrupt the service of a customer in the $M/G/1$ queue. This service is not resumed until all customers that have arrived during the vacation and their offspring have been served (in LCFS order).

We now focus on the amount of work in this $M/G/1$ system at an arbitrary moment \emph{during a visit (busy) period}. Let $K$ be the customer being served at this observation moment, and let $K_A$ be his ancestor.
%Let $K_A$ be the ancestor of the customer being served at an arbitrary moment during a visit period.
By definition, $K_A$ has arrived during a vacation period (or: switch-over period in the corresponding polling system). Denote by $Y_{K_A}$ the amount of work present in the system at the moment that $K_A$ enters the system. An important difference with the situation studied in \cite{boxmagroenendijk87} is that we \emph{cannot} use the PASTA property, so in general $Y_{K_A}\neq Y^{(S)}$. We now condition on the customer type of $K_A$. The mean duration of the service of a type $i$ ancestor and his entire ancestral line is $\E[B_i]/(1-\rho^{(V)})$. %The probability that an arbitrary customer arriving during a switch-over time is a type $i$ customer, is $\frac{\sum_{j=1}^N \lambda_i^{(S_j)} \E[S_j] }{ \sum_{k=1}^N \sum_{j=1}^N \lambda_k^{(S_j)} \E[S_j] }$.
This can be regarded as the mean busy period commencing with the service of an exceptional first customer (namely a type $i$ customer).
Each type $i$ customer arriving during $S_j$, with arrival rate $\lambda_i^{(S_j)}$,  $i,j=1,\dots,N$, starts such a busy period, so the probability that $K_A$ is a type $i$ customer is:
\begin{equation}
p_i = \frac{\sum_{j=1}^N \lambda_i^{(S_j)} \E[S_j] \E[B_i]/(1-\rho^{(V)})}{ \sum_{k=1}^N \sum_{j=1}^N \lambda_k^{(S_j)} \E[S_j] \E[B_k]/(1-\rho^{(V)})}= \frac{\sum_{j=1}^N \lambda_i^{(S_j)} \E[S_j] \E[B_i]}{ \sum_{k=1}^N \sum_{j=1}^N \lambda_k^{(S_j)} \E[S_j] \E[B_k]}.
\label{eq:pitjes}
\end{equation}
Given that $K_A$ is a type $i$ customer, we again pick up the proof of the work decomposition in \cite{boxmagroenendijk87}. Denote by $B_{K_A}$ the service requirement of $K_A$. Then, because of the LCFS service discipline of the $M/G/1$ queue, the amount of work when $K_A$ goes into service is exactly $Y_{K_A}+B_{K_A}$, and the amount of work when the last descendant of $K_A$ has been served equals $Y_{K_A}$ again (for the first time, since the arrival of $K_A$). Ignoring the amount of work present at $K_A$'s arrival, the residual amount of work evolves just as during a busy period in an $M/G/1$ queue with an exceptional first customer (having generic service requirement $B_i$). %The work brought in by the descendants of $K_A$ behaves just like in the visit periods of the polling model. %, e.g. as described in Cases 1 and 2.
The only exception is caused by the vacations (or switch-over times in the polling model), during which the work remains constant or may increase because of new arrivals. However, just as in \cite{boxmagroenendijk87}, if we ignore these vacations and the (LCFS) service of the ancestral lines of the customers that arrive during these vacations, what remains is the workload process during a busy period initiated by a type $i$ customer. Denote by $Y^{(V)}_{M/G/1|i}$ the amount of work at an arbitrary moment during this busy period, and denote by $Y^{(S)}_{A_i}$ the amount of work present in the polling system at an arbitrary \emph{arrival epoch} of a type $i$ customer \emph{during a switch-over time}. Note that $Y_{K_A}$ is distributed like $Y^{(S)}_{A_i}$. Then we have the following decomposition:
\begin{equation}
Y^{(V)} \equaldist Y^{(V)}_{M/G/1|i} + Y^{(S)}_{A_i}  \qquad\text{w.p. }p_i, \qquad i=1,\dots,N,\label{decomp}
\end{equation}
with $p_i$ as given in \eqref{eq:pitjes}, and $Y^{(V)}_{M/G/1|i}$ and $Y^{(S)}_{A_i}$ being independent. This leads to
\begin{equation}
\E[Y^{(V)}] = \sum_{i=1}^N p_i \left(\E[ Y^{(V)}_{M/G/1|i}] + \E[Y^{(S)}_{A_i}]\right),
\label{eq:EVwithoutConditioned}
\end{equation}
with
\begin{align}
\E[Y^{(V)}_{M/G/1|i}] &= \E[R_{B_i}] + \frac{\rho^{(V)}}{1-\rho^{(V)}} \E[R_{B^{(V)}}],\label{eq:EYVmg1i}\\
\E[Y^{(S)}_{A_i}] &= \sum_{j=1}^N \frac{\lambda_i^{(S_j)} \E[S_j]}{\sum_{k=1}^N \lambda_i^{(S_k)} \E[S_k]} \E[Y^{(S_j)}].
\label{eq:EYAi}
\end{align}
%with $\E[Y^{(S_j)}]$ as defined in \eqref{eq:EYSi}. Now we only need to determine $\E[Y^{(V)}_{M/G/1|i}]$.
For \eqref{eq:EYVmg1i} we use standard theory on an $M/G/1$ queue with an exceptional first customer (cf. \cite{wolff1989smt}), and \eqref{eq:EYAi} is established by conditioning on the switch-over period in which a type $i$ customer arrives.

\subsection{PCL for smart customers} % PCL

We are now ready to state the PCL.
\begin{theorem}
%Consider a polling model with smart customers. Under the restriction that the Poisson process according to which work arrives in the system, has a fixed arrival rate during all \emph{visit periods}, and that the amount of work brought by an individual arrival is identically distributed during all \emph{visit periods},
Provided that \eqref{eq:pclassump1} or \eqref{eq:pclassump2} is valid, the following Pseudo-Conservation Law holds:
\begin{align}
\sum_{i=1}^N \overline{\rho}_i \, \E[W_i] &= (1- \overline{\rho}) \sum_{i=1}^N \frac{\E[S_i]}{\E[S]} \E[Y^{(S_i)}] - \sum_{i=1}^N \overline{\rho}_i \E[R_{B_i}] \nonumber\\
&\quad \mbox{}+ \overline{\rho} \sum_{i=1}^N p_i \left(\sum_{j=1}^N \frac{\lambda_i^{(S_j)} \E[S_j]}{\sum_{k=1}^N \lambda_i^{(S_k)} \E[S_k]} \E[Y^{(S_j)}] + \E[R_{B_i}] + \frac{\rho^{(V)}}{1-\rho^{(V)}} \E[R_{B^{(V)}}] \right),
\label{eq:PCLsmart}
\end{align}
where $\E[Y^{(S_i)}]$ are as in~\eqref{eq:EYSi}, and $p_i$ as in~\eqref{eq:pitjes}.

\begin{proof}
We have two equations, \eqref{eq:EVdecomp} and \eqref{EYeqn2}, for the mean total amount of work in the system. Combining these two equations, and plugging in~\eqref{eq:EYinEYSj} and~\eqref{eq:EVwithoutConditioned}, we find
$$\sum_{i=1}^N \left(\E[\LQ_i] \E[B_i] + \overline{\rho}_i \E[R_{B_i}] \right) = (1-\overline{\rho}) \sum_{j=1}^N \frac{\E[S_j]}{\E[S]} \E[Y^{(S_j)}] +
 \overline{\rho} \sum_{i=1}^N p_i \left(\E[ Y^{(V)}_{M/G/1|i}] + \E[Y^{(S)}_{A_i}]\right).$$
By application of Little's law, $\E[\LQ_i] = \overline{\lambda}_i \E[W_i]$, using that $\overline{\rho}_i = \overline{\lambda}_i \E[B_i]$, plugging in \eqref{eq:EYVmg1i} and \eqref{eq:EYAi}, after some rewriting we obtain %the following PCL for a polling model with smart customers: %, given that~\eqref{eq:pclassump1} holds:
%\begin{align}
%\sum_{i=1}^N \overline{\rho}_i \, \E[W_i] &= (1- \overline{\rho}) \sum_{i=1}^N \frac{\E[S_i]}{\E[S]} \E[Y^{(S_i)}] - \sum_{i=1}^N \overline{\rho}_i \E[R_{B_i}] \nonumber\\
%&\quad \mbox{}+ \overline{\rho} \sum_{i=1}^N p_i \left(\sum_{j=1}^N \frac{\lambda_i^{(S_j)} \E[S_j]}{\sum_{k=1}^N \lambda_i^{(S_k)} \E[S_k]} \E[Y^{(S_j)}] + \E[R_{B_i}] + %\frac{\rho^{(V)}}{1-\rho^{(V)}} \E[R_{B^{(V)}}] \right).
%\label{eq:PCLsmart}
%\end{align}
\eqref{eq:PCLsmart}, which is a PCL for a polling model with smart customers.
%
%where $\E[Y^{(S_i)}]$ are as in~\eqref{eq:EYSi}, and the $p_i$ as in~\eqref{eq:pitjes}.
\end{proof}
\end{theorem}

\begin{remark}
When $\lambda_i^{(S_1)} = \lambda_i^{(S_2)} = \ldots = \lambda_i^{(S_N)} = \lambda_i^{(V_1)} = \dots=\lambda_i^{(V_N)} = \lambda_i$, for all $i=1,\ldots,N$, Equation~\eqref{eq:PCLsmart} reduces to~\eqref{eq:pclStandardExh}. E.g., because of PASTA, $\E[Y_{A_i}^{(S)}]=\E[Y^{(S)}]$, and $p_i = \lambda_i/\Lambda$ for all~$i$.
\end{remark}

Case 2, where assumptions \eqref{eq:pclassump2} hold, has a nice practical interpretation if we add the additional requirement that $\sum_{i=1}^N \lambda_i^{(S_j)} = \sum_{i=1}^N \lambda_i^{(V_j)} =: \Lambda$ for all $j=1,\dots,N$. Now, the model can be interpreted as a polling system with customers arriving in one Poisson stream with constant arrival rate $\Lambda$, and generic service requirement $B$, but joining a certain queue with a fixed probability that may depend on the location of the server at the arrival epoch. In Section \ref{examples}, we discuss an example on how these probabilities may be chosen to minimise the mean waiting time of an arbitrary customer. The PCL \eqref{eq:PCLsmart} can be simplified considerably in this situation.
\begin{corollary}
If \eqref{eq:pclassump2} is valid, the PCL \eqref{eq:PCLsmart} reduces to:
\begin{equation}
\sum_{i=1}^N \overline{\rho}_i \, \E[W_i] =  \sum_{i=1}^N \frac{\E[S_i]}{\E[S]} \E[Y^{(S_i)}]  +\frac{\rho^2}{1-\rho}\E[R_B].\label{eq:PCLsmartCaseII}
\end{equation}
\begin{proof}
This is a direct consequence of assumptions \eqref{eq:pclassump2}. E.g., in the computation of \eqref{eq:EYVmg1i} there is no need to condition on a special first customer, and hence the term $\E[Y_{M/G/1|i}]$ does not depend on~$i$ anymore:
\[\E[Y_{M/G/1|i}] = \frac{\E[R_B]}{1-\rho},\]
where $\rho = \Lambda\E[B]$. Additionally, the term $\sum_{i=1}^N p_i \E[Y_{A_i}^{(S)}]$ also simplifies considerably:
\[\sum_{i=1}^N p_i \E[Y_{A_i}^{(S)}] = \sum_{i=1}^N \frac{\E[S_i]}{\E[S]} \E[Y^{(S_i)}].\]
Combining this, multiple terms cancel out and \eqref{eq:PCLsmartCaseII} follows.
It is easily seen that \eqref{eq:PCLsmartCaseII} is in line with~\eqref{eq:pclStandardExh}, when the arrival rates do not change during various visit and switch-over times.
\end{proof}
\end{corollary}

%For the $p_i$ the $\E(B)$ cancel in the nominator and denominator.
%$$
%p_i = \frac{\sum_{j=1}^N \lambda_i^{(S_j)} \E(S_j)}
%           { \sum_{k=1}^N \sum_{j=1}^N \lambda_k^{(S_j)} \E(S_j)}.
%$$
%Combining this gives that the PCL for a polling model with a constant total arrival rate, and equal service time distribution at each of the queues, is given by:
%\begin{equation}
%\sum_{i=1}^N \overline{\rho}_i \, \E[W_i] = (1- \rho) \sum_{i=1}^N \frac{\E[S_i]}{\E(S)} \E(Y_{S_i})  +\frac{\rho}{1-\rho}\E(R_B)
%+ \rho\, \frac{\sum_{i=1}^N \sum_{j=1}^N \lambda_i^{(S_j)} \E(S_j) \E(Y_{S_j})}
%           { \sum_{k=1}^N \sum_{l=1}^N \lambda_k^{(S_l)} \E[S_l]}.
%\label{eq:PCLsmartCaseII}
%\end{equation}
%where we have used that in this case $\sum_{i=1}^N \overline{\rho}_i \E[R_{B_i}] = \overline{\rho} \, \E(R_B)$, and in the last term, $\sum_{j=1}^N \lambda_i^{(S_j)} \E(S_j)$ of the $p_i$ cancelled out against such a term in front of the $\E(Y_{S_j})$.

%[checked for a few numerical examples with 3 queues, in all of the cases turned out to be correct]

\section{Numerical examples\label{examples}}

\subsection{Example 1: smart customers}

\begin{figure}[h!]
\begin{center}
\includegraphics[width=0.67\textwidth]{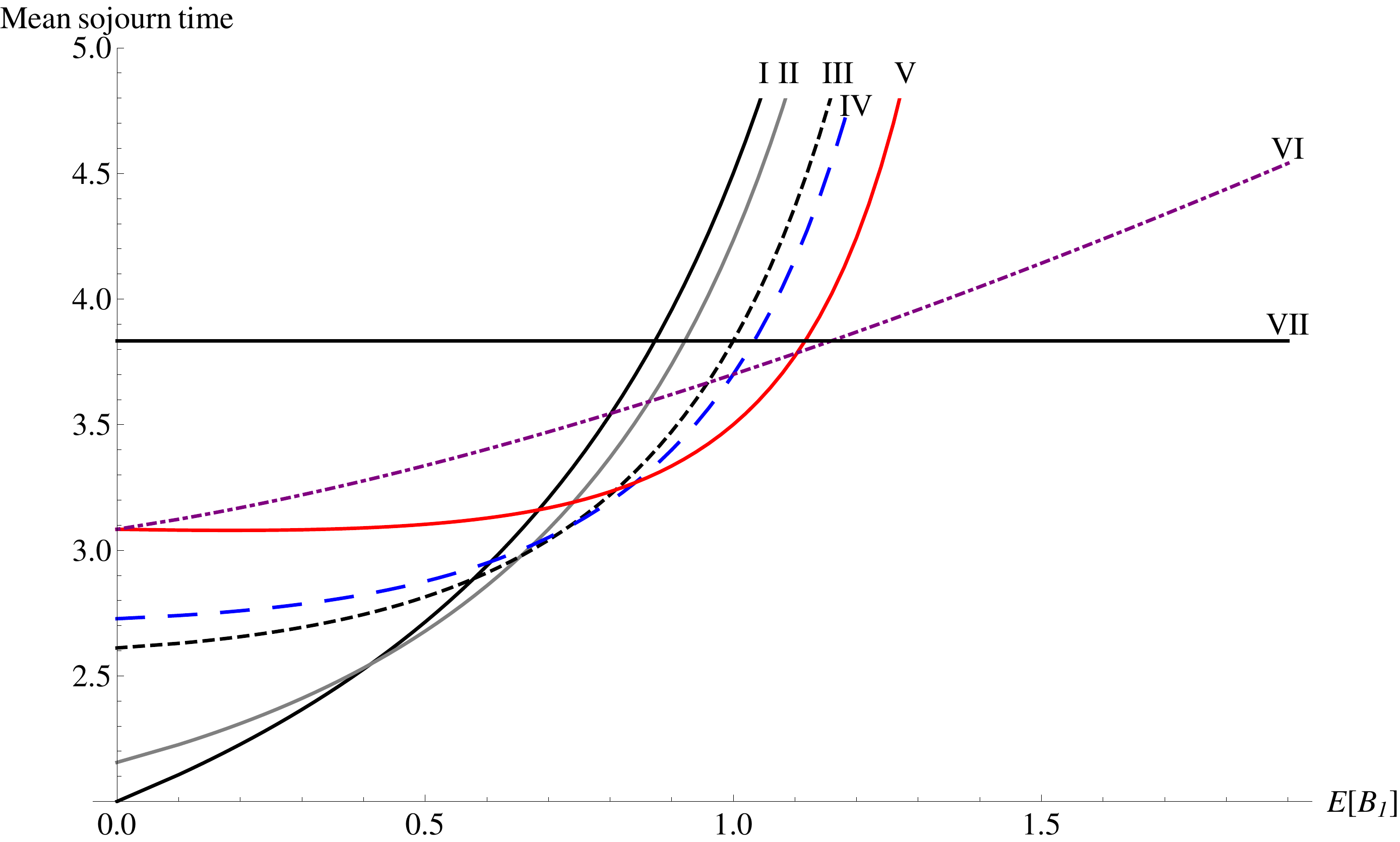}
\end{center}
\caption{The mean sojourn time of an arbitrary customer for the seven smartest strategies in Example 1, against the mean service time in $Q_1$.}
\label{figsmartstrategies}
\end{figure}

In the first numerical example, we study a polling system where arriving customers choose which queue they join, based on the current position of the server. In \cite{smartcustomers,boxmakelbert94} a fully symmetric case is studied with gated service, and it is proven that the mean sojourn time of customers is minimised if customers join the queue that is being served directly after the queue that is currently being served. Although the exhaustive case is not studied, it is intuitively clear that in this situation smart customers join the queue that is currently being served. Or, in case an arrival takes place during a switch-over time, join the next queue that is visited. In this example, we study this situation in more detail by adding an extra parameter that can be varied. The polling model is fully symmetric, except for the service time of customers in $Q_1$, which is varied. The practical interpretation is the following: as in the previously described examples, customers arrive with a fixed arrival intensity, say $\Lambda$, and choose which queue they join. This does not affect their service time, except when they choose $Q_1$. In this case the service time has a different distribution. To illustrate the dynamics of this system, we choose the following setting. The system consists of three queues with exhaustive service. The switch-over times are all exponentially distributed with mean $1$. The service times are also exponentially distributed with $\E[B_2] = \E[B_3] = 1$, and $\E[B_1]$ is varied between 0 and 2. Arriving customers choose one queue which they want to join. This queue is the same for all customers, so there is no randomness involved in the selection, which is only based on the location of the server at their arrival epochs. We intend to find the optimal queue for customers to join. In terms of the model parameters: we seek to find values for $\lambda_i^{(V_j)}$ and $\lambda_i^{(S_j)}$, $i,j=1,2,3$, that minimise the mean sojourn time of an arbitrary customer, under the restriction that for each value of $j$, exactly one $\lambda_i^{(V_j)}$ and exactly one $\lambda_i^{(S_j)}$ is equal to $\Lambda$, and all the other values are 0. A valid combination of these arrival intensities is called a \emph{strategy}, and we introduce the short notation for a strategy by the indices of the queues that are joined in respectively $(V_1, S_1, V_2, S_2, V_3, S_3)$. E.g., for the fully symmetric case, with $\E[B_1] = 1$, it is intuitively clear that the optimal strategy is to join $Q_{i}$, if the arrival takes place during $V_i$, and to join $Q_{i+1}$ if the arrival takes place during $S_i$. This strategy is denoted by $(1, 2, 2, 3, 3, 1)$, and corresponds to $\lambda_1^{(V_1)}=\lambda_2^{(V_2)}=\lambda_3^{(V_3)}=\Lambda$, and $\lambda_2^{(S_1)}=\lambda_3^{(S_2)}=\lambda_1^{(S_3)}=\Lambda$. The other arrival intensities are 0.
As stated before, we vary $\E[B_1]$ between 0 and 2, and focus on the overall mean sojourn time. It is clear that making $\E[B_1]$ smaller, makes it more attractive to join $Q_1$ (even if another queue is served), whereas making $\E[B_1]$ larger, makes it less attractive to join $Q_1$. In order to obtain numerical results, we choose the (arbitrary) value $\Lambda=\frac35$. It turns out that as much as \emph{seven} different strategies can be optimal, depending on the value of $\E[B_1]$. We refer to these strategies as I through VII, listed in Table \ref{smartstrategies}, along with their region of optimality. For each of these strategies, the mean sojourn time of an arbitrary customer is plotted versus $\E[B_1]$ in Figure \ref{figsmartstrategies}.

\begin{table}[h!]
\begin{center}
\begin{tabular}{|c|cccccc|c|}
\hline
Strategy & \multicolumn{6}{|c|}{Queue to join during} & Region of optimality \\
& $V_1$ & $S_1$ & $V_2$ & $S_2$ & $V_3$ & $S_3$ & \\
\hline
I &   1 & 1 & X & 1 & X & 1 & $0.00    \leq \E[B_1] \leq 0.41$ \\
II &  1 & 2 & 1 & 1 & X & 1 & $0.41 \leq \E[B_1] \leq 0.66$ \\
III & 1 & 2 & 2 & 1 & X & 1 & $0.66 \leq \E[B_1] \leq 0.73$\\
IV &  1 & 2 & 2 & 3 & 1 & 1 & $0.73 \leq \E[B_1] \leq 0.84$\\
V &   1 & 2 & 2 & 3 & 3 & 1 & $0.84 \leq \E[B_1] \leq 1.10$\\
VI &  2 & 2 & 2 & 3 & 3 & 1 & $1.10 \leq \E[B_1] \leq 1.16$\\
VII & X & 2 & 2 & 3 & 3 & 2 & $1.16 \leq \E[B_1] \phantom{\leq 0.00\,} $\\
\hline
\end{tabular}
\end{center}
\caption{The seven smartest strategies in Example 1 that minimise the mean waiting time of an arbitrary customer who can choose the queue in which he wants to be served. An `X' means that the length of the corresponding visit time equals 0 because customers never join this queue.}
\label{smartstrategies}
\end{table}

As expected, $Q_1$ is most popular if $\E[B_1]$ is very small. In particular, for very small values of $\E[B_1]$, customers \emph{always} join this queue (Strategy~I). As $\E[B_1]$ becomes larger, $Q_2$ and later also $Q_3$ are chosen in more and more situations (Strategies II--V). Strategy~V, which is optimal if the system is (nearly) symmetric, is the one where customers join the queue that is being served, or is going to be served next if the arrival takes place during a switch-over time. Strategy~VI, which is optimal in only a very small range of values of $\E[B_1]$, states that customers only join $Q_1$ during the switch-over time $S_3$. Strategy~VII, in which customers never join $Q_1$, is optimal for large values of $\E[B_1]$. The ergodicity constraint, considering all parameters are fixed except for $\E[B_1]$, for the different strategies is also interesting to mention. For strategies I-V, the necessary and sufficient condition for stability is $\E[B_1] < \frac53$. Strategies~VI and~VII always result in a stable system, regardless of $\E[B_1]$. For illustration purposes, we show how to compute the ergodicity constraint for Strategy~V. As indicated in Section \ref{nosubtypes}, the ergodicity constraint requires computation of the eigenvalues of the matrix $R - I_N$, where $I_N$ is the $N \times N$ identity matrix, and $R$ is an $N \times N$ matrix containing elements $\rho_{ij} := \lambda_{i}^{(V_j)}\E[B_i]$. For Strategy~V, we find
\[
R-I_3 = \begin{pmatrix}
\Lambda\E[B_1] - 1& 0 & 0 \\
0 & \Lambda-1 & 0 \\
0 & 0 & \Lambda-1
\end{pmatrix}.
\]
The eigenvalues of this matrix are $\Lambda\E[B_1] - 1, \Lambda-1$, and again $\Lambda-1$. The ergodicity constraint in this situation states that the largest (and, hence, all) of these eigenvalues should be negative. This means that $\E[B_1]<\frac{1}{\Lambda}$ is a sufficient and necessary condition for stability of this system, given that $\Lambda<1$. The ergodicity constraints of the other strategies are computed similarly, but all rows and columns corresponding to visit times that are zero should be deleted from the matrix $R-I_3$ (cf. \cite{pollinglevy09}). For Strategies {VI and VII} this implies that the first row and the first column should be deleted. Since the first row is the only row which contains $\E[B_1]$, these strategies always result in a stable system. Note that the arrival rates during switch-over times do not play a role in the ergodicity constraint.

It is also interesting to discuss what \emph{stupid customers} would do in this system. Stupid customers choose the worst possible strategy, in order to maximise the mean sojourn time of an arbitrary customer. We do not go into details and do not mention exactly which strategy is worst for each value of $\E[B_1]$, but we pick out some interesting cases. Obviously, when $\E[B_1] = 0$, the worst possible thing to do is never to join $Q_1$. The worst strategy in this case is $(X, 3, 3,2,2,3)$, where $X$ means that any queue can be chosen (because the length of the corresponding visit time equals 0, since customers never join this queue). This strategy leads to an overall mean sojourn time of $7.48$.
As $\E[B_1]$ grows larger, $Q_1$ gradually will be chosen more frequently. In the symmetric case, $\E[B_1] = 1$, customers arriving during $V_i$ choose $Q_{i-1}$, and customers arriving during $S_i$ choose $Q_i$, resulting in a mean sojourn time of $8.5$. For large $\E[B_1]$, the worst possible strategy might be a bit surprising. It is \emph{not} simply to always join $Q_1$, but it is $(1,1,1,2,1,3)$. During visit periods, customers always join $Q_1$, but during $S_i$ customers join $Q_i$. For $\E[B_1] \uparrow \frac53$, this strategy results in the highest mean sojourn time of an arbitrary customer. %The fact that $\E[B_1]$ is higher than the mean service time in the other two queues is not enough reason  for customers arriving during switch-over periods to join $Q_1$. This is caused by the fact that the difference between the mean service time in $Q_1$ and the mean service time in the other two queues is not that big, because $\E[B_1] < \frac53$.
For the situation $\E[B_1] \geq \frac53$, there are many strategies for which the system becomes unstable and sojourn times become infinite. The worst possible strategy for $\E[B_1] \geq \frac53$ that still results in a stable system, is $(3, 1, X, 1, 1, 1)$.

\subsection{Example 2: no arrivals during a specific period}

In this example we illustrate how to deal with polling models with arrival rates being zero during certain periods. For MVA, this is no problem. The equations presented in Section \ref{mva} still give the correct solution if some of the arrival rates during periods are zero. The problem arises when determining the LST of the waiting time distribution \eqref{waitingtimeLST} and can only be circumvented by a work-around, which is explained using a simple example. The polling model in this example contains two queues, $Q_1$ and $Q_2$, which are served exhaustively. All switch-over times and all service times are exponentially distributed with parameter 1. All arrival rates are $\frac12$, except for the arrival rate of type 1 customers arriving during the service of type 2 customers: $\lambda_1^{(V_2)} = 0$. This brings along some complications. First of all, \eqref{cycletimeLSTexh} cannot be used to determine the cycle time LST. This is no real problem, because \eqref{cycletimeLST} can be used instead. Because of $\lambda_1^{(V_2)}$ being zero, we should use \eqref{queuelengthGFduringvisitjLambda0} instead of \eqref{queuelengthGFduringvisitj} for type $1^{(V_2)}$ customers to determine the PGF of the steady-state queue length of $Q_1$. Again, no real problem but just something to be careful about. Determining the waiting time LST for type 1 customers does raise some issues, though. The (generalisation of the) distributional form of Little's law, given by \eqref{waitingtimeLST}, uses the joint distribution of customers left behind by a departing type $i$ customer to determine his time spent in the system. As can be seen in the proof of Theorem \ref{distformLittleTheorem}, this technique requires that type~$i$ customers may arrive during each period within a cycle. In our model this is not the case, because no type 1 customers arrive during $V_2$. This implies that the number of customers left behind by a departing type $1$ customer, does not give any information about the waiting time of type $1$ customers (more specifically, of those that arrived during $S_1$), because a departing type 1 customer does not leave behind any customers (of any type) that have arrived during $V_2$.

A work-around for this problem, is to introduce an \emph{extra queue}, $Q_{X}$, with type $X$ customers that have no service requirement $(B_X=0)$, and $\lambda_X^{(V_2)}>0$. Customers in this queue are served exhaustively somewhere between the end of $V_1$ and the beginning of $V_2$, because type $X^{(V_2)}$ customers have to be present at departure epochs of type $1$ customers. In our approach, we choose to treat $Q_X$ as a regular queue between $Q_1$ and $Q_2$ with no switch-over time from $Q_X$ to $Q_2$ because this gives us a ``normal'', three-queue polling system. Determining the waiting time LST of type 1 customers, requires a careful application of the distributional form of Little's law to the various customer subtypes in Equation \eqref{jointQLatdeparture}. For convenience, we introduce the following two vectors, where the elements correspond to customer subtypes $(1^{(V_1)},\dots,1^{(S_2)}, X^{(V_1)},\dots,X^{(S_2)}, 2^{(V_1)},\dots,2^{(S_2)})$:
\begin{align*}
%\bm{\omega_1} &= (1-\frac{\omega}{\lambda_1^{(V_1)}}, 1-\frac{\omega}{\lambda_1^{(S_1)}}, 1, 1-\frac{\omega}{\lambda_1^{(S_2)}}, 1, \dots, 1),\\
%\bm{\omega^*_1} &= (1-\frac{\omega}{\lambda_1^{(V_1)}}, 1-\frac{\omega}{\lambda_1^{(S_1)}}, 1, 1-\frac{\omega}{\lambda_1^{(S_2)}}, 1, \dots, 1,1-\frac{\omega}{\lambda_X^{(V_2)}},1,\dots, 1),
\bm{\omega_1} &= (1-\frac{\omega}{\lambda_1^{(V_1)}}, 1-\frac{\omega}{\lambda_1^{(S_1)}}, 1, 1-\frac{\omega}{\lambda_1^{(S_2)}}, 1, 1, 1, 1, 1, 1, 1, 1),\\
\bm{\omega^*_1} &= (1-\frac{\omega}{\lambda_1^{(V_1)}}, 1-\frac{\omega}{\lambda_1^{(S_1)}}, 1, 1-\frac{\omega}{\lambda_1^{(S_2)}}, 1, 1,1-\frac{\omega}{\lambda_X^{(V_2)}},1,1, 1, 1, 1),
\end{align*}
the difference being in the element corresponding to the type $X$ customers that arrive during $V_2$. Note that we do not introduce customer subtypes that arrive during $V_X$ or $S_X$, because the lengths of these periods are 0. The LST of the waiting time distribution of type 1 customers is given by:
\[
\E\left[\ee^{-\omega W_1}\right] = \frac{1}{\beta_1(\omega)} \frac{\overline{\lambda}}{\overline{\lambda}_1}\left(M_1^{(V_1)}(\bm{\omega_1})+M_1^{(S_1)}(\bm{\omega^*_1})+M_1^{(V_2)}(\bm{\omega_1})+M_1^{(S_2)}(\bm{\omega_1})\right).
\]
The interpretation is that we use the type $X^{(V_2)}$ customers left behind by a departing $1^{(S_1)}$ customer to determine the length of $V_2$, which is part of the total waiting time of a type $1^{(S_1)}$ customer. The other type 1 customers arrive after the visit to $Q_2$ and can be handled in the regular way. The numerical results of this example are shown in Table \ref{example2numericalResults}.
\begin{table}[h!]
\begin{center}
\begin{tabular}{|l|r|r|}
\hline
& $Q_1$ & $Q_2$ \\
\hline
Mean queue length at arrival epochs & 1.750  & 3.375  \\
Mean queue length at departure epochs &  1.750 & 3.375  \\
Mean queue length at arbitrary epochs &  1.188 & 3.375  \\
Mean waiting time &  3.750 & 5.750  \\
Standard deviation waiting time & 5.093  & 6.280  \\
\hline
\end{tabular}
\end{center}
\caption{Numerical results for the polling model discussed in Example 2.}
\label{example2numericalResults}
\end{table}

We can modify \eqref{cycletimeLSTexh} and \eqref{intervisittimeLSTexh} accordingly to obtain the LSTs of the cycle time distribution $C_1^*$, starting at a visit ending to $Q_1$, and the intervisit time distribution $I_1$:
\begin{align*}
%\E\big[\ee^{-\omega C^*_1}\big] &= V_{b_1}^{(S_1)}\big(\pi_1(\omega)-\frac{\omega}{\lambda_1^{(V_1)}}, \pi_1(\omega)-\frac{\omega}{\lambda_1^{(S_1)}}, 1, \pi_1(\omega)-\frac{\omega}{\lambda_1^{(S_2)}}, 1, \dots, 1,1-\frac{\omega}{\lambda_X^{(V_2)}},1,\dots, 1\big),\\
\E\big[\ee^{-\omega C^*_1}\big] &= \VB_{1}^{(S_1)}\big(\pi_1(\omega)-\frac{\omega}{\lambda_1^{(V_1)}}, \pi_1(\omega)-\frac{\omega}{\lambda_1^{(S_1)}}, 1, \pi_1(\omega)-\frac{\omega}{\lambda_1^{(S_2)}}, 1, 1,1-\frac{\omega}{\lambda_X^{(V_2)}},1,1,1,1, 1\big),\\
\E\big[\ee^{-\omega I_1}\big] &= \VB_{1}^{(S_1)}\big(\bm{\omega^*_1}\big).\\
\end{align*}

\appendix
\section*{Appendix}
\section{MVA equations}\label{appendixMVAeqns}

In this appendix we present all MVA equations that have been omitted in Section \ref{mva}.
%eqns4b

The mean duration of the next period~$V_i$, when in~$S_j$ is denoted by $\EintForw{V_i}{S_j}$. A difference with $\EintForw{V_i}{V_j}$, is that $\EintForw{V_i}{S_i}$ is not different from $\EintForw{V_i}{S_j}$ for $j\neq i$. Similar to \eqref{eqmva_EViVjf}, we have for $i = 1,\dots,N$, $j = i, \ldots, i+N-1$:
\begin{align}
\EintForw{V_i}{S_j} &= \E[\BP_i] \left(\E[\LQ_i^{(S_j)}] + \lambda_i^{(S_j)} \E[R_{S_j}]
                    + \sum_{k=j+1}^{i+N-1} \left(\lambda_i^{(V_k)} \EintForw{V_k}{S_j} +
                    \lambda_i^{(S_k)} \E[S_k]  \right)\right).
\label{eqmva_EViSjf}\\
\intertext{%
Equation \eqref{eqmva_ERSiVj_forw} for $\E[R_{S_i : V_j}]$, the mean residual duration of the interval $S_i, V_{i+1}, \ldots, V_j$, is obtained by conditioning on the period in which the interval is observed, looking forwards in time. Similarly, we find expressions for $\E[R_{S_i : S_j} ], \E[R_{V_i : V_j}]$, and $\E[R_{V_i : S_j} ]$. For $i = 1,\dots,N$, $j = i+1,\dots,i+N-1$:}
%eqns33b
\E[R_{S_i : S_j} ] &= \sum_{k=i}^{j} \frac{\E[S_k]}{\E[(S_i\!:\!S_j)]}\left( \E[R_{S_k}] + \sum_{l=k+1}^j\left( \E[S_l] + \EintForw{V_l}{S_k} \right)\right)\nonumber\\
                    &\quad \mbox{} + \sum_{k=i+1}^{j} \frac{\E[V_k]}{\E[(S_i\!:\!S_j)]}\left( \sum_{l=k}^{j}\left( \E[S_l] + \EintForw{V_l}{V_k} \right)\right).
\label{eqmva_ERSiSj_forw}\\
\intertext{For $i = 1,\dots,N$, $j = i+1,\dots,i+N-1$:}
%eqns33c
\E[R_{V_i : V_j}] &= \sum_{k=i}^{j-1} \frac{\E[S_k]}{\E[(V_i\!:\!V_j)]}\left(  \E[R_{S_k}] + \EintForw{V_j}{S_k} + \sum_{l=k+1}^{j-1}\left( \E[S_l] +  \EintForw{V_l}{S_k} \right)\right)\nonumber\\
                    &\quad \mbox{} + \sum_{k=i}^{j} \frac{\E[V_k]}{\E[(V_i\!:\!V_j)]}\left( \EintForw{V_j}{V_k} +  \sum_{l=k}^{j-1}\left( \E[S_l] + \EintForw{V_l}{V_k} \right)\right).
\label{eqmva_ERViVj_forw}\\
\intertext{For $i = 1,\dots,N$, $j = i+1,\dots,i+N-1$:}
%eqns33d
\E[R_{V_i : S_j} ] &= \sum_{k=i}^j \frac{\E[S_k]}{\E[(V_i\!:\!S_j)]}\left( \E[R_{S_k}] + \sum_{l=k+1}^j\left( \E[S_l] + \EintForw{V_l}{S_k} \right)\right)\nonumber\\
                    &\quad \mbox{} + \sum_{k=i}^j \frac{\E[V_k]}{\E[(V_i\!:\!S_j)]}\left( \sum_{l=k}^j \left(\E[S_l] + \EintForw{V_l}{V_k} \right)\right).
\label{eqmva_ERViSj_forw}
\\
\intertext{%
In Section \ref{mva}, a second set of equations is discussed for $\E[R_{S_i : V_j} ], \E[R_{S_i : S_j} ], \E[R_{V_i : V_j}]$, and $\E[R_{V_i : S_j} ]$. This set is obtained by conditioning on the period in which the interval is observed, but now looking backwards in time. We use that the residual length of an interval has the same distribution as the elapsed time of this interval. The equation for $\E[R_{S_i : V_j} ]$ is given by \eqref{eqmva_ERSiVj}. The other equations are given below. For $i = 1,\dots,N$, $j = i+1,\dots,i+N-1$:}
%eqns3b
\E[R_{S_i : S_j} ] &= \sum_{k=i}^{j} \frac{\E[S_k]}{\E[(S_i\!:\!S_j)]}\left( \EintBack{S_i}{S_k} + \sum_{l=i+1}^k \left(\EintBack{S_l}{S_k} + \EintBack{V_l}{S_k} \right)\right)\nonumber\\
                    &\quad \mbox{} + \sum_{k=i+1}^{j} \frac{\E[V_k]}{\E[(S_i\!:\!S_j)]}\left( \EintBack{S_i}{V_k} + \EintBack{V_k}{V_k} + \sum_{l=i+1}^{k-1} \left(\EintBack{S_l}{V_k} + \EintBack{V_l}{V_k} \right)\right).
\label{eqmva_ERSiSj}
\intertext{For $i = 1,\dots,N$, $j = i+1,\dots,i+N-1$:}
%eqns3c
\E[R_{V_i : V_j} ] &= \sum_{k=i}^{j-1} \frac{\E[S_k]}{\E[(V_i\!:\!V_j)]}\left(  \sum_{l=i}^k \left( \EintBack{S_l}{S_k} + \EintBack{V_l}{S_k} \right)\right)\nonumber\\
                    &\quad \mbox{} + \sum_{k=i}^{j} \frac{\E[V_k]}{\E[(V_i\!:\!V_j)]}\left( \EintBack{V_k}{V_k} +  \sum_{l=i}^{k-1} \left(\EintBack{S_l}{V_k} + \EintBack{V_l}{V_k} \right)\right).
\label{eqmva_ERViVj}\\
\intertext{For $i = 1,\dots,N$, $j = i,\dots,i+N-1$:}
%eqns3d
\E[R_{V_i : S_j} ] &= \sum_{k=i}^{j} \frac{\E[S_k]}{\E[(V_i\!:\!S_j)]}\left(  \sum_{l=i}^k \left( \EintBack{S_l}{S_k} + \EintBack{V_l}{S_k} \right)\right)\nonumber\\
                    &\quad \mbox{} + \sum_{k=i}^{j} \frac{\E[V_k]}{\E[(V_i\!:\!S_j)]}\left( \EintBack{V_k}{V_k} +  \sum_{l=i}^{k-1} \left( \EintBack{S_l}{V_k} + \EintBack{V_l}{V_k} \right)\right).
\label{eqmva_ERViSj}
\end{align}

\expandafter\ifx\csname urlstyle\endcsname\relax
  \providecommand{\doi}[1]{DOI: #1}\else
  \providecommand{\doi}{DOI: \begingroup \urlstyle{rm}\Url}\fi

%\bibliographystyle{abbrvnat}
%\bibliography{smartcustRevision}

\end{document}